\documentclass[8pt,a4wide]{article}

\usepackage[normalem]{ulem}
\usepackage{listings}
\usepackage{showlabels}
\usepackage{a4wide}
\usepackage{amsmath} 
\usepackage{amsthm}
\usepackage{amsfonts}
\usepackage{amssymb,amscd}
\usepackage{multirow}
\usepackage[mathcal]{eucal}
\usepackage[spanish,USenglish]{babel} 
\usepackage{graphics}  
\usepackage{graphicx} 
\usepackage{hyperref}
\usepackage{fancyhdr}
\usepackage{titlesec}
\usepackage[latin1]{inputenc}
\usepackage{amsmath}
\usepackage{amsfonts}
\usepackage{amssymb}
\usepackage{amsthm}
\usepackage{newlfont}
\usepackage{amsbsy}
\usepackage{bm}
\usepackage[all]{xy}
\usepackage{textcomp}
\usepackage{color}
\usepackage{epsfig}
\usepackage{dsfont}
\usepackage{latexsym}
\usepackage{array}
\usepackage{longtable}
\usepackage{enumerate}
\usepackage{multicol}
\usepackage{mathrsfs}
\usepackage{cancel}
\usepackage{wasysym}
\theoremstyle{plain}
\newtheorem{teor}{Theorem}
\newtheorem{lemma}{Lemma}
\newtheorem{coro}{Corollary}
\newtheorem{prop}{Proposition}

\theoremstyle{definition}

\newtheorem*{ejem*}{Examples}
\newtheorem{defin}{Definition}
\newtheorem*{demos}{Proof}
\newtheorem{remark}{Remark}
\theoremstyle{remark}

\newcommand{\complex}{\mathbf{\mathbb{C}}}

\newcommand{\ente}{\mathbf{\mathbb{Z}}}
\newcommand{\natu}{\mathbf{\mathbb{N}}}
\newcommand{\llave}[1]{\left\{ #1\right\}}

\newcommand{\corch}[1]{\left[ #1\right]}
\newcommand{\paren}[1]{\left( #1\right)}

\usepackage[breakable]{tcolorbox}
\usepackage{parskip} 

\usepackage{iftex}
\ifPDFTeX
\usepackage[T1]{fontenc}
\usepackage{mathpazo}
\else
\usepackage{fontspec}
\fi

\usepackage{graphicx}

\usepackage{caption}
\DeclareCaptionFormat{nocaption}{}
\usepackage[Export]{adjustbox} 
\adjustboxset{max size={0.9\linewidth}{0.9\paperheight}}
\usepackage{float}
\usepackage{xcolor} 
\usepackage{enumerate} 
\usepackage{geometry} 
\usepackage{amsmath} 
\usepackage{amssymb} 
\usepackage{textcomp} 
\AtBeginDocument{%
}
\usepackage{upquote} 
\usepackage{eurosym} 
\usepackage[mathletters]{ucs} 
\usepackage{fancyvrb} 
\usepackage{grffile} 
\makeatletter 
\def\Gread@@xetex#1{%
	\IfFileExists{"\Gin@base".bb}%
	{\Gread@eps{\Gin@base.bb}}%
	{\Gread@@xetex@aux#1}%
}
\makeatother

\usepackage{hyperref}
\usepackage{titling}
\usepackage{longtable} 
\usepackage{booktabs}  
\usepackage[inline]{enumitem} 
\usepackage[normalem]{ulem} 
\usepackage{mathrsfs}

\definecolor{urlcolor}{rgb}{0,.145,.698}
\definecolor{linkcolor}{rgb}{.71,0.21,0.01}
\definecolor{citecolor}{rgb}{.12,.54,.11}

\definecolor{ansi-black}{HTML}{3E424D}
\definecolor{ansi-black-intense}{HTML}{282C36}
\definecolor{ansi-red}{HTML}{E75C58}
\definecolor{ansi-red-intense}{HTML}{B22B31}
\definecolor{ansi-green}{HTML}{00A250}
\definecolor{ansi-green-intense}{HTML}{007427}
\definecolor{ansi-yellow}{HTML}{DDB62B}
\definecolor{ansi-yellow-intense}{HTML}{B27D12}
\definecolor{ansi-blue}{HTML}{208FFB}
\definecolor{ansi-blue-intense}{HTML}{0065CA}
\definecolor{ansi-magenta}{HTML}{D160C4}
\definecolor{ansi-magenta-intense}{HTML}{A03196}
\definecolor{ansi-cyan}{HTML}{60C6C8}
\definecolor{ansi-cyan-intense}{HTML}{258F8F}
\definecolor{ansi-white}{HTML}{C5C1B4}
\definecolor{ansi-white-intense}{HTML}{A1A6B2}
\definecolor{ansi-default-inverse-fg}{HTML}{FFFFFF}
\definecolor{ansi-default-inverse-bg}{HTML}{000000}

\let\Oldtex\TeX
\let\Oldlatex\LaTeX
\renewcommand{\TeX}{\textrm{\Oldtex}}
\renewcommand{\LaTeX}{\textrm{\Oldlatex}}
\title{First Example}

\makeatletter
\def\PY@reset{\let\PY@it=\relax \let\PY@bf=\relax%
	\let\PY@ul=\relax \let\PY@tc=\relax%
	\let\PY@bc=\relax \let\PY@ff=\relax}
\def\PY@tok#1{\csname PY@tok@#1\endcsname}
\def\PY@toks#1+{\ifx\relax#1\empty\else%
	\PY@tok{#1}\expandafter\PY@toks\fi}
\def\PY@do#1{\PY@bc{\PY@tc{\PY@ul{%
				\PY@it{\PY@bf{\PY@ff{#1}}}}}}}
\def\PY#1#2{\PY@reset\PY@toks#1+\relax+\PY@do{#2}}

\expandafter\def\csname PY@tok@w\endcsname{\def\PY@tc##1{\textcolor[rgb]{0.73,0.73,0.73}{##1}}}
\expandafter\def\csname PY@tok@c\endcsname{\let\PY@it=\textit\def\PY@tc##1{\textcolor[rgb]{0.25,0.50,0.50}{##1}}}
\expandafter\def\csname PY@tok@cp\endcsname{\def\PY@tc##1{\textcolor[rgb]{0.74,0.48,0.00}{##1}}}
\expandafter\def\csname PY@tok@k\endcsname{\let\PY@bf=\textbf\def\PY@tc##1{\textcolor[rgb]{0.00,0.50,0.00}{##1}}}
\expandafter\def\csname PY@tok@kp\endcsname{\def\PY@tc##1{\textcolor[rgb]{0.00,0.50,0.00}{##1}}}
\expandafter\def\csname PY@tok@kt\endcsname{\def\PY@tc##1{\textcolor[rgb]{0.69,0.00,0.25}{##1}}}
\expandafter\def\csname PY@tok@o\endcsname{\def\PY@tc##1{\textcolor[rgb]{0.40,0.40,0.40}{##1}}}
\expandafter\def\csname PY@tok@ow\endcsname{\let\PY@bf=\textbf\def\PY@tc##1{\textcolor[rgb]{0.67,0.13,1.00}{##1}}}
\expandafter\def\csname PY@tok@nb\endcsname{\def\PY@tc##1{\textcolor[rgb]{0.00,0.50,0.00}{##1}}}
\expandafter\def\csname PY@tok@nf\endcsname{\def\PY@tc##1{\textcolor[rgb]{0.00,0.00,1.00}{##1}}}
\expandafter\def\csname PY@tok@nc\endcsname{\let\PY@bf=\textbf\def\PY@tc##1{\textcolor[rgb]{0.00,0.00,1.00}{##1}}}
\expandafter\def\csname PY@tok@nn\endcsname{\let\PY@bf=\textbf\def\PY@tc##1{\textcolor[rgb]{0.00,0.00,1.00}{##1}}}
\expandafter\def\csname PY@tok@ne\endcsname{\let\PY@bf=\textbf\def\PY@tc##1{\textcolor[rgb]{0.82,0.25,0.23}{##1}}}
\expandafter\def\csname PY@tok@nv\endcsname{\def\PY@tc##1{\textcolor[rgb]{0.10,0.09,0.49}{##1}}}
\expandafter\def\csname PY@tok@no\endcsname{\def\PY@tc##1{\textcolor[rgb]{0.53,0.00,0.00}{##1}}}
\expandafter\def\csname PY@tok@nl\endcsname{\def\PY@tc##1{\textcolor[rgb]{0.63,0.63,0.00}{##1}}}
\expandafter\def\csname PY@tok@ni\endcsname{\let\PY@bf=\textbf\def\PY@tc##1{\textcolor[rgb]{0.60,0.60,0.60}{##1}}}
\expandafter\def\csname PY@tok@na\endcsname{\def\PY@tc##1{\textcolor[rgb]{0.49,0.56,0.16}{##1}}}
\expandafter\def\csname PY@tok@nt\endcsname{\let\PY@bf=\textbf\def\PY@tc##1{\textcolor[rgb]{0.00,0.50,0.00}{##1}}}
\expandafter\def\csname PY@tok@nd\endcsname{\def\PY@tc##1{\textcolor[rgb]{0.67,0.13,1.00}{##1}}}
\expandafter\def\csname PY@tok@s\endcsname{\def\PY@tc##1{\textcolor[rgb]{0.73,0.13,0.13}{##1}}}
\expandafter\def\csname PY@tok@sd\endcsname{\let\PY@it=\textit\def\PY@tc##1{\textcolor[rgb]{0.73,0.13,0.13}{##1}}}
\expandafter\def\csname PY@tok@si\endcsname{\let\PY@bf=\textbf\def\PY@tc##1{\textcolor[rgb]{0.73,0.40,0.53}{##1}}}
\expandafter\def\csname PY@tok@se\endcsname{\let\PY@bf=\textbf\def\PY@tc##1{\textcolor[rgb]{0.73,0.40,0.13}{##1}}}
\expandafter\def\csname PY@tok@sr\endcsname{\def\PY@tc##1{\textcolor[rgb]{0.73,0.40,0.53}{##1}}}
\expandafter\def\csname PY@tok@ss\endcsname{\def\PY@tc##1{\textcolor[rgb]{0.10,0.09,0.49}{##1}}}
\expandafter\def\csname PY@tok@sx\endcsname{\def\PY@tc##1{\textcolor[rgb]{0.00,0.50,0.00}{##1}}}
\expandafter\def\csname PY@tok@m\endcsname{\def\PY@tc##1{\textcolor[rgb]{0.40,0.40,0.40}{##1}}}
\expandafter\def\csname PY@tok@gh\endcsname{\let\PY@bf=\textbf\def\PY@tc##1{\textcolor[rgb]{0.00,0.00,0.50}{##1}}}
\expandafter\def\csname PY@tok@gu\endcsname{\let\PY@bf=\textbf\def\PY@tc##1{\textcolor[rgb]{0.50,0.00,0.50}{##1}}}
\expandafter\def\csname PY@tok@gd\endcsname{\def\PY@tc##1{\textcolor[rgb]{0.63,0.00,0.00}{##1}}}
\expandafter\def\csname PY@tok@gi\endcsname{\def\PY@tc##1{\textcolor[rgb]{0.00,0.63,0.00}{##1}}}
\expandafter\def\csname PY@tok@gr\endcsname{\def\PY@tc##1{\textcolor[rgb]{1.00,0.00,0.00}{##1}}}
\expandafter\def\csname PY@tok@ge\endcsname{\let\PY@it=\textit}
\expandafter\def\csname PY@tok@gs\endcsname{\let\PY@bf=\textbf}
\expandafter\def\csname PY@tok@gp\endcsname{\let\PY@bf=\textbf\def\PY@tc##1{\textcolor[rgb]{0.00,0.00,0.50}{##1}}}
\expandafter\def\csname PY@tok@go\endcsname{\def\PY@tc##1{\textcolor[rgb]{0.53,0.53,0.53}{##1}}}
\expandafter\def\csname PY@tok@gt\endcsname{\def\PY@tc##1{\textcolor[rgb]{0.00,0.27,0.87}{##1}}}
\expandafter\def\csname PY@tok@err\endcsname{\def\PY@bc##1{\setlength{\fboxsep}{0pt}\fcolorbox[rgb]{1.00,0.00,0.00}{1,1,1}{\strut ##1}}}
\expandafter\def\csname PY@tok@kc\endcsname{\let\PY@bf=\textbf\def\PY@tc##1{\textcolor[rgb]{0.00,0.50,0.00}{##1}}}
\expandafter\def\csname PY@tok@kd\endcsname{\let\PY@bf=\textbf\def\PY@tc##1{\textcolor[rgb]{0.00,0.50,0.00}{##1}}}
\expandafter\def\csname PY@tok@kn\endcsname{\let\PY@bf=\textbf\def\PY@tc##1{\textcolor[rgb]{0.00,0.50,0.00}{##1}}}
\expandafter\def\csname PY@tok@kr\endcsname{\let\PY@bf=\textbf\def\PY@tc##1{\textcolor[rgb]{0.00,0.50,0.00}{##1}}}
\expandafter\def\csname PY@tok@bp\endcsname{\def\PY@tc##1{\textcolor[rgb]{0.00,0.50,0.00}{##1}}}
\expandafter\def\csname PY@tok@fm\endcsname{\def\PY@tc##1{\textcolor[rgb]{0.00,0.00,1.00}{##1}}}
\expandafter\def\csname PY@tok@vc\endcsname{\def\PY@tc##1{\textcolor[rgb]{0.10,0.09,0.49}{##1}}}
\expandafter\def\csname PY@tok@vg\endcsname{\def\PY@tc##1{\textcolor[rgb]{0.10,0.09,0.49}{##1}}}
\expandafter\def\csname PY@tok@vi\endcsname{\def\PY@tc##1{\textcolor[rgb]{0.10,0.09,0.49}{##1}}}
\expandafter\def\csname PY@tok@vm\endcsname{\def\PY@tc##1{\textcolor[rgb]{0.10,0.09,0.49}{##1}}}
\expandafter\def\csname PY@tok@sa\endcsname{\def\PY@tc##1{\textcolor[rgb]{0.73,0.13,0.13}{##1}}}
\expandafter\def\csname PY@tok@sb\endcsname{\def\PY@tc##1{\textcolor[rgb]{0.73,0.13,0.13}{##1}}}
\expandafter\def\csname PY@tok@sc\endcsname{\def\PY@tc##1{\textcolor[rgb]{0.73,0.13,0.13}{##1}}}
\expandafter\def\csname PY@tok@dl\endcsname{\def\PY@tc##1{\textcolor[rgb]{0.73,0.13,0.13}{##1}}}
\expandafter\def\csname PY@tok@s2\endcsname{\def\PY@tc##1{\textcolor[rgb]{0.73,0.13,0.13}{##1}}}
\expandafter\def\csname PY@tok@sh\endcsname{\def\PY@tc##1{\textcolor[rgb]{0.73,0.13,0.13}{##1}}}
\expandafter\def\csname PY@tok@s1\endcsname{\def\PY@tc##1{\textcolor[rgb]{0.73,0.13,0.13}{##1}}}
\expandafter\def\csname PY@tok@mb\endcsname{\def\PY@tc##1{\textcolor[rgb]{0.40,0.40,0.40}{##1}}}
\expandafter\def\csname PY@tok@mf\endcsname{\def\PY@tc##1{\textcolor[rgb]{0.40,0.40,0.40}{##1}}}
\expandafter\def\csname PY@tok@mh\endcsname{\def\PY@tc##1{\textcolor[rgb]{0.40,0.40,0.40}{##1}}}
\expandafter\def\csname PY@tok@mi\endcsname{\def\PY@tc##1{\textcolor[rgb]{0.40,0.40,0.40}{##1}}}
\expandafter\def\csname PY@tok@il\endcsname{\def\PY@tc##1{\textcolor[rgb]{0.40,0.40,0.40}{##1}}}
\expandafter\def\csname PY@tok@mo\endcsname{\def\PY@tc##1{\textcolor[rgb]{0.40,0.40,0.40}{##1}}}
\expandafter\def\csname PY@tok@ch\endcsname{\let\PY@it=\textit\def\PY@tc##1{\textcolor[rgb]{0.25,0.50,0.50}{##1}}}
\expandafter\def\csname PY@tok@cm\endcsname{\let\PY@it=\textit\def\PY@tc##1{\textcolor[rgb]{0.25,0.50,0.50}{##1}}}
\expandafter\def\csname PY@tok@cpf\endcsname{\let\PY@it=\textit\def\PY@tc##1{\textcolor[rgb]{0.25,0.50,0.50}{##1}}}
\expandafter\def\csname PY@tok@c1\endcsname{\let\PY@it=\textit\def\PY@tc##1{\textcolor[rgb]{0.25,0.50,0.50}{##1}}}
\expandafter\def\csname PY@tok@cs\endcsname{\let\PY@it=\textit\def\PY@tc##1{\textcolor[rgb]{0.25,0.50,0.50}{##1}}}

\makeatother

\makeatletter
\newbox\Wrappedcontinuationbox 
\newbox\Wrappedvisiblespacebox 
\newcommand*\Wrappedvisiblespace {\textcolor{red}{\textvisiblespace}} 
\newcommand*\Wrappedcontinuationsymbol {\textcolor{red}{\llap{\tiny$\m@th\hookrightarrow$}}} 
\newcommand*\Wrappedcontinuationindent {3ex } 
\newcommand*\Wrappedafterbreak {\kern\Wrappedcontinuationindent\copy\Wrappedcontinuationbox} 
\newcommand*\Wrappedbreaksatspecials {%
	\def\PYGZus{\discretionary{\char`\_}{\Wrappedafterbreak}{\char`\_}}%
	\def\PYGZob{\discretionary{}{\Wrappedafterbreak\char`\{}{\char`\{}}%
	\def\PYGZcb{\discretionary{\char`\}}{\Wrappedafterbreak}{\char`\}}}%
	\def\PYGZca{\discretionary{\char`\^}{\Wrappedafterbreak}{\char`\^}}%
	\def\PYGZam{\discretionary{\char`\&}{\Wrappedafterbreak}{\char`\&}}%
	\def\PYGZlt{\discretionary{}{\Wrappedafterbreak\char`\<}{\char`\<}}%
	\def\PYGZgt{\discretionary{\char`\>}{\Wrappedafterbreak}{\char`\>}}%
	\def\PYGZsh{\discretionary{}{\Wrappedafterbreak\char`\#}{\char`\#}}%
	\def\PYGZpc{\discretionary{}{\Wrappedafterbreak\char`\%}{\char`\%}}%
	\def\PYGZdl{\discretionary{}{\Wrappedafterbreak\char`\$}{\char`\$}}%
	\def\PYGZhy{\discretionary{\char`\-}{\Wrappedafterbreak}{\char`\-}}%
	\def\PYGZsq{\discretionary{}{\Wrappedafterbreak\textquotesingle}{\textquotesingle}}%
	\def\PYGZdq{\discretionary{}{\Wrappedafterbreak\char`\"}{\char`\"}}%
	\def\PYGZti{\discretionary{\char`\~}{\Wrappedafterbreak}{\char`\~}}%
} 
\newcommand*\Wrappedbreaksatpunct {%
	\lccode`\~`\.\lowercase{\def~}{\discretionary{\hbox{\char`\.}}{\Wrappedafterbreak}{\hbox{\char`\.}}}%
	\lccode`\~`\,\lowercase{\def~}{\discretionary{\hbox{\char`\,}}{\Wrappedafterbreak}{\hbox{\char`\,}}}%
	\lccode`\~`\;\lowercase{\def~}{\discretionary{\hbox{\char`\;}}{\Wrappedafterbreak}{\hbox{\char`\;}}}%
	\lccode`\~`\:\lowercase{\def~}{\discretionary{\hbox{\char`\:}}{\Wrappedafterbreak}{\hbox{\char`\:}}}%
	\lccode`\~`\?\lowercase{\def~}{\discretionary{\hbox{\char`\?}}{\Wrappedafterbreak}{\hbox{\char`\?}}}%
	\lccode`\~`\!\lowercase{\def~}{\discretionary{\hbox{\char`\!}}{\Wrappedafterbreak}{\hbox{\char`\!}}}%
	\lccode`\~`\/\lowercase{\def~}{\discretionary{\hbox{\char`\/}}{\Wrappedafterbreak}{\hbox{\char`\/}}}%
	\catcode`\.\active
	\catcode`\,\active 
	\catcode`\;\active
	\catcode`\:\active
	\catcode`\?\active
	\catcode`\!\active
	\catcode`\/\active 
	\lccode`\~`\~ 	
}
\makeatother

\let\OriginalVerbatim=\Verbatim
\makeatletter
\renewcommand{\Verbatim}[1][1]{%
	\sbox\Wrappedcontinuationbox {\Wrappedcontinuationsymbol}%
	\sbox\Wrappedvisiblespacebox {\FV@SetupFont\Wrappedvisiblespace}%
	\def\FancyVerbFormatLine ##1{\hsize\linewidth
		\vtop{\raggedright\hyphenpenalty\z@\exhyphenpenalty\z@
			\doublehyphendemerits\z@\finalhyphendemerits\z@
			\strut ##1\strut}%
	}%
	\def\FV@Space {%
		\nobreak\hskip\z@ plus\fontdimen3\font minus\fontdimen4\font
		\discretionary{\copy\Wrappedvisiblespacebox}{\Wrappedafterbreak}
		{\kern\fontdimen2\font}%
	}%
	
	\Wrappedbreaksatspecials
	\OriginalVerbatim[#1,codes*=\Wrappedbreaksatpunct]%
}
\makeatother

\definecolor{incolor}{HTML}{303F9F}
\definecolor{outcolor}{HTML}{D84315}
\definecolor{cellborder}{HTML}{CFCFCF}
\definecolor{cellbackground}{HTML}{F7F7F7}

\makeatletter
\newcommand{\boxspacing}{\kern\kvtcb@left@rule\kern\kvtcb@boxsep}
\makeatother

\hypersetup{
	breaklinks=true,  
	colorlinks=true,
	urlcolor=urlcolor,
	linkcolor=linkcolor,
	citecolor=citecolor,
}

\begin{document}
	\title{\textbf{Matrix Bispectrality of Full Rank One Algebras}}
	\author{Brian D. Vasquez C. & Jorge P. Zubelli}
	\author{
Brian D. Vasquez$^{1}$, Jorge P. Zubelli$^{2}$\\ \\
\small{$^{1}$IMPA, $^{2}$Department of Mathematics, Khalifa University.}\\
\small{\texttt{$^{1}$bridava927@gmail.com, $^{2}$zubelli@gmail.com}}
}
\date{\small{\today}}
	\maketitle

\begin{abstract}
We study algebraic properties of full rank $ 1 $ algebras in a general framework and derive a method to verify if one such matrix polynomial  sub-algebra is bispectral. We give two examples illustrating the method. In the first one, we consider the eigenvalue to be scalar-valued, whereas in the second one,  we assume it to be matrix-valued. In the former example we put forth a Pierce decomposition of that algebra. \\

As a byproduct, we answer positively a conjecture of F.~A.~Gr\"{u}nbaum concerning certain  noncommutative matrix algebras associated to the bispectral problem.\\

\end{abstract}

\linespread{2}

\section{Introduction}

Classical orthogonal polynomials as many important special functions satisfy remarkable relations both in the physical as well as in the spectral variables~\cite{DuranGrunbaum2004}. More precisely, they are eigenfunctions of an operator in the physical variable (say $x$) with eigenvalues depending on the spectral variable (say $z$) as well as the other way around, eigenfunctions of an operator in $z$ with  
$x$--dependent eigenvalues. 
Such bispectral property was explored in the scalar case in the work of J.~J.~Duistermaat and F.~A.~Gr\"{u}nbaum \cite{duistermaat1986differential}. It turned out to have deep connections with many problems in Mathematical Physics. Indeed, it could be 
arranged in suitable manifolds which were naturally parameterized by the flows of the Korteweg de-Vries (KdV) hierarchy or its master-symmetries
\cite{zubelli1991differential,Zubelli2000}. It led to generalizations associated to the 
Kadomtsev-Petviashvili (KP) hierarchy~\cite{Zubelli1992,wilson,horozov2002bispectral,iliev1999discrete}.

%
%
%
%
This article concerns a promising generalization of the original bispectral problem. More precisely, we consider the triples $(L,\psi, B)$ satisfying systems of equations 
	\begin{equation}\label{bispec}
	L\psi(x,z)=\psi (x,z)F(z) \hspace{1 cm}  (\psi B) (x,z)=\theta(x)\psi(x,z)
	\end{equation}
	with $L=L(x,\partial_{x})$, $B=B(z,\partial_{z})$  linear matrix  differential operators, i.e., $L\psi=\sum_{i=0}^{l}a_{i}(x) \cdot \partial_{x}^i \psi$,
	$\psi B=\sum_{j=0}^{m}\partial_{z}^j \psi \cdot b_{j}(z)$. The functions $a_{i}, b_{j},F, \theta$ and the nontrivial common eigenfunction $\psi$ are in principle compatible sized matrix valued functions.
	
A triple $(L,\psi, B)$ satisfying \eqref{bispec} is called a {\em bispectral triple}. We remark that all the differential operators are considered in  a neighborhood of an arbitrary given point and following \cite{duistermaat1986differential} we assume that the functions are smooth enough so that all the derivatives considered make sense. 

Now we fix the normalized~\footnote{ If $L=L(x,\partial_{x})$, $L=\sum_{i=0}^{l}a_{i}(x)\partial_{x}^i$ with $a_{l}$ constant and scalar, $a_{l-1}=0$, then $L$ is called normalized.} operator $L$ and the eigenfunctions $\psi(\cdot,z)$. We are interested in the bispectral pairs associated to $L=L(x,\partial_{x})$, i.e., operators  $B=B(z, \partial_{z})$ such that $(\psi B)(x,z)=\theta(x)\psi(x,z)$ for some function $\theta=\theta(x)$. It is not hard to verify that the set of operators  $B=B(z,\partial_{z})$ satisfying \eqref{bispec} generates a noncommutative
algebra of operators.

The main goal of this article is to establish a method to verify whether an algebra of matrix polynomials is bispectral or not. For a given scalar eigenvalue function the corresponding algebra of  matrix eigenvalues is characterized. Moreover, the isomorphism between the matrix eigenvalues and the corresponding operator is given explicitly. These results may be used to give positive answers to the three conjectures in~\cite{Grfrm-4}.

The noncommutative (or matrix) version of the bispectral problem was first studied in \cite{zubelli1990differential,Zubelli1992a,Zubelli1992b} for the situation where  both the physical and  spectral operators were acting on the same side of the eigenfunction and the eigenvalues are scalar valued. Later on, several generalizations were considered. 
See~\cite{SZ01, Kas15, BL08, GI03, Grfrm-4, GHY17} and 
references therein. This article follows up on the possibility of having the 
physical and the spectral operators acting on different sides.  
We also follow the suggestion in \cite{Grfrm-4} of considering both eigenvalues as matrix valued.

The plan of this article is as follows: In Section~\ref{sec1}, we develop some preliminary theoretical results which allow us to determine if an algebra is bispectral and a subalgebra of the matrix polynomial algebra. In Section~\ref{sec2} we apply these results to the three algebras presented in \cite{Grfrm-4} to prove that they are bispectral. In  Subsection~\ref{first} the function $F(z)$  is scalar valued and $\theta(x)$ is matrix valued, while in the Subsection~\ref{spin} both are matrix valued. 

 In the Appendix, we prove the noncommutative Ad-Condition which implies that the bispectral algebras are subalgebras of the matrix polynomial algebra.
 Furthermore, some computations were performed with the software Singular and Maxima.

\section{General Results}\label{sec1}

We consider the triples $(L,\psi, B)$ satisfying systems of equations 
	\begin{equation}\label{bispec}
	L\psi(x,z)=\psi (x,z)F(z) \hspace{1 cm}  (\psi B) (x,z)=\theta(x)\psi(x,z)
	\end{equation}
	with $L=L(x,\partial_{x})$, $B=B(z,\partial_{z})$  linear matrix  differential operators, i.e., $L\psi=\sum_{i=0}^{l}a_{i}(x) \cdot \partial_{x}^i \psi$,
	$\psi B=\sum_{j=0}^{m}\partial_{z}^j \psi \cdot b_{j}(z)$. The functions $a_{i}:U \subset \complex \rightarrow M_{N}(\complex) , b_{j}:V\subset \complex \rightarrow M_{N}(\complex),F:V\subset \complex \rightarrow M_{N}(\complex), \theta:U\subset \complex\rightarrow M_{N}(\complex)$ and the nontrivial common eigenfunction 
	$\psi:U\times V\subset \complex^2\rightarrow M_{N}(\complex)$ are in principle compatible sized meromorphic matrix valued functions defined in suitable open subsets $U,V\subset \complex$. 
	
A triple $(L,\psi, B)$ satisfying \eqref{bispec} is called a bispectral triple.

Now we fix the normalized operator $L$ and the eigenfunctions $\psi(\cdot,z)$. We are interested in the bispectral pairs associated to $L=L(x,\partial_{x})$, i.e., operators  $B=B(z, \partial_{z})$ such that $(\psi B)(x,z)=\theta(x)\psi(x,z)$ for some function $\theta=\theta(x)$. It is not hard to verify that the set of operators  $B=B(z,\partial_{z})$ satisfying \eqref{bispec} generates a noncommutative
algebra of operators.

We first note that $\theta$ satisfying Equation \eqref{bispec} has to be an element of  the algebra of polynomials with $N\times N$ matrix coefficients, which we denote by $M_{N}(\complex)\corch{x}$. The proof follows closely an argument in the original paper of \cite{duistermaat1986differential}. 
See the Appendix~\ref{adcond}.\\
Clearly the set 
\begin{equation}
\mathbb{A}=\llave{\theta \in  M_{N}(\complex)\corch{x} \big| \exists \mathcal{B}=\mathcal{B}(z,\partial_{z}), (\psi\mathcal{B})(x,z)=\theta(x)\psi(x,z)}
\end{equation}
 is a  noncommutative $\complex$-algebra.

We shall start with some theoretical results. 	

\begin{defin}
	Define $Bp(z,L)$ to be the set of 
	{\em bispectral partners} 
	to $L$, i.e.,
\begin{equation*}
Bp(z,L)=\llave{B=B(z,\partial_{z}) | \exists \theta\in M_{N}(\complex[x]), (\psi\mathcal{B})(x,z)=\theta(x)\psi(x,z)} \mbox{ .}
\end{equation*}

\end{defin} 

A straightforward consequence of the definition is the following.

\begin{lemma}
	The set $Bp(z,L)$ is a $\complex$-algebra.
\end{lemma}

However, much more can be said about the properties of the algebra $Bp(z,L)$ in the case that will be studied in the sequel. For that we have to consider the following important class of algebras:
\begin{defin}
	Let $\mathbb{K}$ be a field, $C$ be a graded $\mathbb{K}$-algebra, we define a full rank $1$ algebra to be a subalgebra $A\subset C$ such that
	\begin{equation*}
	A=E\oplus \bigoplus_{j=k_{0}}^{\infty} C_{j}
	\end{equation*} 
for some finite dimensional $\mathbb{K}$-vector space $E$ and $k_{0}\in \natu$. Furthermore, we denote by $\llave{e_{i}^{[j]}}_{1\leq i\leq N, j\geq 0}$ some basis for $C_j$. See \cite{Casper2017}.
\end{defin}

\begin{remark}
Note that for a full rank 1 algebra $A\subset C$. 
We consider $k_{0}$ the smallest positive integer such that $ C_{j}\subset A$, for all $j\geq k_{0}$. For this $k_0$ we can write $E=
\paren{\bigoplus_{j=0}^{k_{0}-1} C_{j}}\cap A$.
\end{remark}

The results in Theorems \ref{full}, \ref{Noether} and \ref{bis} will be used in the sequel to provide a positive answer to the  conjectures of Grunbaum \cite{Grfrm-4}. They are of interest on their own.

\begin{teor}\label{full}
	Let $C$ be a graded $\mathbb{K}$-algebra where  $\dim_{\mathbb{K}} C_{j}=N< \infty$, $C_{j}=\sum_{i=1}^{N} \mathbb{K}\cdot e_{i}^{[j]} $. Suppose that for every $t$, $1\leq t\leq N$, $ i,j\in \natu$, there exist $1\leq r,s\leq N$ such that $e_{t}^{[i+j]}=e_{r}^{[i]}e_{s}^{[j]}$ and $A\subset C$ is a full rank $1$ algebra. Then, $A$ is a finitely generated $\mathbb{K}$-algebra. 
\end{teor}
\begin{demos}
We write \begin{equation*}
A=E\oplus \bigoplus_{j=k_{0}}^{\infty} C_{j}  \mbox{\rm  .}
\end{equation*}
Since $E$ is a finite dimensional $\mathbb{K}$-vector space $E$, we can consider a basis $\llave{\alpha_{1},...,\alpha_{m}}$ for $E$ and write 
$E=\sum_{s=1}^{m}\mathbb{K}\cdot \alpha_{s}$. Define 
\begin{equation*}
A_{0}:=\mathbb{K}\cdot\langle e_{i}^{[k]},\alpha_{s}\mid 1\leq i \leq N, k_{0}\leq k \leq 2k_{0}-1, 1\leq s \leq m\rangle
\mbox{\rm  .}
\end{equation*}
 We claim that $A=A_{0}$.
 
 First of all, we prove that for every $q\geq 2$, $(q-1)k_{0}\leq k \leq qk_{0}-1$, $e_{i}^{[k]}\in A_{0}$. The initial step is clear for $q=2$. Assume that $e_{i}^{[p]}\in A_{0}$ for $(q-1)k_{0}\leq p \leq qk_{0}-1$, $1\leq i \leq N$ and note that $qk_{0}\leq k \leq (q+1)k_{0}-1$ implies $(q-1)k_{0}\leq k-k_{0}\leq qk_{0}-1$ and $e_{i}^{[k-k_0]}\in A_{0}$ for $1\leq i \leq N$. Consider $1\leq i \leq N$, by hypothesis there exists $1\leq r,s\leq N$ such that $e_{i}^{[k]}=e_{r}^{[k-k_{0}]}\cdot e_{s}^{k_{0}} \in A_{0}$. This proves the inductive step. The assertion follows by induction.

Since $k_{0}+\natu=\bigcup_{q=2}^{\infty}\llave{k\in \natu \mid (q-1)k_{0}\leq k \leq qk_{0}-1}$. We have that $e_{i}^{[k]}\in A_{0}$ for $1\leq i \leq N$,
$k\geq k_{0}$ then $\bigoplus_{j=k_{0}}^{\infty} C_{j}\subset A_{0}$. But $E=\sum_{s=1}^{m}\mathbb{K}\cdot \alpha_{s}\subset A_{0}$. Thus, $A=A_{0}$ and $A$ is a finitely generated $k$-algebra.
\end{demos}
\begin{remark}
The converse is not true. Consider for example the graded algebra $C=M_{N}(\mathbb{K}[x])$ and $A=\mathbb{K}[x]$, then $A$ is a finitely generated $\mathbb{K}$-algebra which is not of full rank 1.
 
\end{remark}
Now we use the following theorem whose proof may be found in \cite{Stafford}.
\begin{teor}[Stafford] \label{Noether}
	Let $R\subset S$ be algebras over a central field $\mathbb{K}$ such that $S$ is Noetherian and $S/R$ is a finite dimensional $\mathbb{K}$-vector space. Then, $R$ is Noetherian.
\end{teor}
\begin{coro}
	Let $C$ be a Noetherian graded algebra and $A\subset C$ be a full rank 1 $\mathbb{K}$-algebra over a central field $\mathbb{K}$. Then, $A$ is Noetherian.
\end{coro}
\begin{demos}
	Since $A=E\oplus \bigoplus_{j=k_{0}}^{\infty} C_{j}$ for some finite dimensional vector space $E$, we can consider the complement of any subspace $F$ with respect to $\bigoplus_{j=0}^{k_{0}-1} C_{j}$ and obtain $C=F\oplus A$ then $\dim_{\mathbb{K}}(C/A)=\dim_{\mathbb{K}} (F)< \infty$. Since $C$ is Noetherian the previous theorem implies the assertion. \qed
\end{demos}

\begin{defin}
Let  $\mathbb{K}$ be a field, $C$ be a $\mathbb{K}$-algebra and $S\subset C$. We define 
$$\mathbb{K}\cdot <S>=span\llave{\prod_{j=1}^{n}s_j\mid s_1, ... , s_n \in S, n \in \natu } \mbox{ .}$$
\end{defin} 

The following theorem connects  the bispectral property to full rank $1$ algebras.

\begin{teor}[Full Rank One Algebras]\label{bis}
	Let $C$ be a graded $\mathbb{K}$-algebra, $\Gamma \subset C$ and $A\subset C$ two $\mathbb{K}$-algebras with the following properties:
	
	\begin{enumerate}
		\item 
		$\Gamma$ is a full rank $1$ algebra, with decomposition $\Gamma=E\oplus \bigoplus_{j=k_{0}}^{\infty}C_{j}$, for some $k_{0}\in \natu$. 
		\item  $\mathbb{K}\cdot \langle E\rangle=\Gamma$.
		\item  $A\bigcap \oplus_{k=0}^{k_{0}-1}C_{k}=E$.
	\end{enumerate}
	Then, $\Gamma=A$. 
\end{teor} 
\begin{demos} We shall break the proof in 2 steps.
	
	Step 1: The inclusion $\Gamma\subset A$.\\
	Using (2) and (3) we have $E\subset A$ and $\Gamma$ is the algebra generated by $E$, since $A$ is an algebra we obtain the inclusion  $\Gamma\subset A$.
	
	Step 2: The inclusion $A \subset \Gamma$.\\
	Consider $\theta \in A$ and write $\theta=\theta_{1}+\theta_2$ with $\theta_{1}\in \oplus_{k=0}^{k_{0}-1}C_{k}$ and  $\theta_{2}\in \oplus_{k=k_{0}}^{\infty}C_{k}$, since $\Gamma \supset  \oplus_{k=k_{0}}^{\infty}C_{k}$ we have that $\theta_{2}\in \Gamma\subset A$. In particular $\theta_{1}= \theta -\theta_{2}\in A\bigcap \oplus_{k=0}^{k_{0}-1}C_{k}=E\subset \Gamma$, then $\theta=\theta_{1}+\theta_2 \in \Gamma$. $\qed$
\end{demos}

\begin{defin}
	The shift operator $S_{N}\in M_{N}(\mathbb{K}[x])$ is defined by  $$S_{N}=\sum_{s=1}^{N-1}e_{s,s+1}$$ 
	for $N\geq 2$, where as usual
	$e_{r,s}$ denotes the matrix with $1$ at entry $(r,s)$ and zeros elsewhere.  
\end{defin}
We recall that
for $N \geq 2$ 
	\begin{equation*}
	S_{N}^j = \left\{
	\begin{array}{ll}
	\sum_{s=1}^{N-j} e_{s,s+j} & \mbox{\rm if\ } 0\leq j \leq N-1, \\
	0 & \mbox{\rm if\ } j \geq N.
	\end{array}
	\right.
	\end{equation*}
	In particular $S_{N}$ is nilpotent of degree $N$.
	
	The following theorem give us a concrete example of a nontrivial full rank $1$ algebra.
\begin{teor}\label{example}
	Let $N\in \ente_{+}$ and the following elements in  $M_{N}(\mathbb{K}[x])$:\\
$	\alpha_{0}=S_{N}$, \\
$\alpha_{1}=Ix+(-1)^{N}e_{N1}x^{N}$,\\
	$\alpha_{k}=e_{1N}x^{k}, \mbox{\rm if\ } 2\leq k\leq N-1$, \\
		$\beta_{k}=e_{kk}x^N+(-1)^N e_{N1}x^{2N-1}, \mbox{\rm if\ } 1\leq k\leq N$.
		
Then, $x^{2N}M_{N}(\mathbb{K}[x])$ is contained in the subalgebra $ A$ of $M_{N}(\mathbb{K}[x])$ that is generated by $\alpha_{j},\beta_{k}$, $0\leq j \leq N-1$, $1\leq k \leq N$.
\end{teor}
\begin{demos}
Since $\beta_{k}^2=e_{kk}x^{2N}$ for $2\leq k \leq N-1$ and $\alpha_{1}^{n}=Ix^{n}+(-1)^{n}ne_{N1}x^{n+N-1}$ for $n\geq 1$ we have
$\beta_{k}\alpha_{1}^{n}=e_{kk}x^{n+2N}\in A$ for $n\geq 1$, in other words $e_{kk}x^{n}\in A$ for $n\geq 2N$.
On the other hand, $\alpha_{0}^{N-1}\alpha_{1}^{n}\alpha_{0}^{N-1}=(-1)^{N}ne_{1N}x^{n+N-1}\in A$ for $n\geq 1$, hence $e_{1N}x^{n}\in A$ for $n\geq N$. However, $e_{1N}x^k\in A$ for $2\leq k \leq N-1$, therefore $e_{1N}x^{n}\in A$ for $n\geq 2$. 

Note that $\alpha_{2}\alpha_{1}^n=e_{1N}x^n+(-1)^N n e_{11}x^{n+N-1}\in A$, $\alpha_{1}^n \alpha_{2}=e_{1N}x^n+(-1)^N n e_{NN}x^{n+N-1}\in A$. Then $e_{11}x^n \in 	A$ and 
$e_{NN}x^n \in 	A$ for $n \geq N+1$.
This implies that $\beta_{1}\alpha_{1}^n=(e_{11}x^N+(-1)^N e_{N1}x^{2N-1})(Ix^n+(-1)^N n e_{N1}x^{n+N-1})=
e_{11}x^{n+N}+(-1)^N e_{N1}x^{2N+n-1}\in A$ for $n\geq 1$. Thus,  $e_{N1}x^n \in A$ for $n\geq 2N$.

The previous proposition implies $\alpha_{0}^j =\sum_{s=1}^{N-j}e_{s,s+j}$ for $0\leq j\leq N-1$ then $\alpha_{0}^{N-i}(e_{N1}x^n)\alpha_{0}^{j-1}=e_{ij}x^n \in A$ for $1\leq 
i,j \leq N$, $n\geq 2N$ and this proves the assertion.\qed
 \end{demos}

\begin{coro}
The algebra $A$ is full rank $1$
\end{coro}
\begin{demos}
	Note that $A=E\oplus x^{2N}M_{N}(\mathbb{K}[x])$ with $E=A\cap \oplus_{j=0}^{2N-1}M_{N}(\mathbb{K}[x])_{j}$. \qed
\end{demos}

In the next section we give a family of algebras whose bispectrality can be obtained using the Theorem \ref{bis}.

\section{The Examples} \label{sec2}
We begin with the example of the matrix algebra given in the paper \cite{Grfrm-4}. There the algebra considered is the set of polynomials of the form 
\begin{equation}\label{motivation}
\theta(x)=\paren{
	\begin{matrix} 
	r_{0}^{11} & r_{0}^{12} \\
	0 & r_{0}^{11} 
	\end{matrix}}+  \paren{
	\begin{matrix} 
	r_{1}^{11} & r_{1}^{12} \\
	0 & r_{1}^{11} 
	\end{matrix}}x+ \paren{
	\begin{matrix} 
	r_{2}^{11} & r_{2}^{12} \\
	r_{1}^{11} & r_{2}^{22} 
	\end{matrix}}x^{2}+
\paren{
	\begin{matrix} 
	r_{3}^{11} & r_{3}^{12} \\
	r_{2}^{22}+r_{2}^{11}-r_{1}^{12} & r_{3}^{22} 
	\end{matrix}}x^3+x^{4}p(x)
\end{equation} 
where $p\in M_{2}(\complex)[x]$ and all the variables $r_{0}^{11},r_{0}^{12}, r_{1}^{11},r_{1}^{12},r_{2}^{11},r_{2}^{22},r_{3}^{11},r_{3}^{12}, r_{3}^{22} \in \complex$. Note that this algebra is full rank one and the relations that must be determined to obtain the complete description of the algebra are in the monomials of degree less than or equal to three. In the following subsection, we generalize this algebra to an arbitrary size of matrix $N$ and find the relations that determine them.

\subsection{Family of Algebras Linked to a Nilpotent Element in $M_{N}(\mathbb{K})$ }\label{first}

As a particular example, we consider a nilpotent element $S\in M_{N}(\mathbb{K})$ of degree $D\geq 2$, consider the matrix valued function 
\begin{equation*}
\psi(x,z)=
e^{xz}\paren{Iz+\sum_{m=1}^{D}(-1)^{m} S^{m-1}x^{-m}} \mbox{ , }
\end{equation*}
and note that $L\psi(x,z)=-z^2\psi(x,z)$ for the ordinary differential operator
\begin{equation*}
L= -\partial_{x}^2+2
\sum_{m=1}^{D}(-1)^{m+1} mS^{m-1}x^{-m-1}.
\end{equation*}
 
Moreover, if $m$ is even we have:
\begin{equation*}
(ad L)^{m}(\theta)\psi=(-1)^{m/2}2^{m}(-z^{2})^{m/2}\psi b_{m}
=(-1)^{m/2}2^{m}L^{m/2}\psi b_{m}= \paren{(-1)^{m/2}2^{m}(L^{m/2})\cdot b_{m}}\psi.
\end{equation*}
Therefore,  
\begin{equation*}
\paren{(ad L)^{m}(\theta)- \paren{(-1)^{m/2}2^{m}(L^{m/2})\cdot b_{m}}}\psi=0.
\end{equation*}
However, the operator $(ad L)^{m}(\theta)- \paren{(-1)^{m/2}2^{m}(L^{m/2})\cdot b_{m}}$ is independent of $z$ and its kernel contains the infinite dimensional linearly independent set   $\llave{\psi(\cdot,z)}_{z\in \complex}$. Thus,  the operator is zero and $(ad L)^{m}(\theta)= (-1)^{m/2}2^{m}(L^{m/2})\cdot b_{m}$.

Now we characterize the algebra $\mathbb{A}=\llave{\theta \in M_{N}(\mathbb{K}[x])| \exists \mathcal{B}=\mathcal{B}(z,\partial_{z}), (\psi\mathcal{B})(x,z)=\theta(x)\psi(x,z)}$ for this particular example. We begin with the definition of the family $\mathcal{P}=\llave{P_{k}}_{k\in \natu}$ which will be used to describe  the map $\theta \mapsto \mathcal{B}$ such that 
$(\psi\mathcal{B})(x,z)=\theta(x)\psi(x,z)$.

\begin{defin}
	For $k\in \natu$ and $\theta\in M_{N}(\mathbb{K}[x])$, we define 
	\begin{equation}
	P_{k}(\theta)=\frac{\theta^{(k+1)}(0)}{k!}-\sum_{j=k+2}^{k+D}(-1)^{k-j}\corch{\frac{\theta^{(j)}(0)}{j!},S^{j-k-1}}.
	\end{equation}
\end{defin}

The family $\mathcal{P}=\llave{P_k}_{k\in \natu}$ can be used to describe the algebra prescribed by Equation~\eqref{motivation}. If $\theta(x)=\sum_{k=0}^{m}a_k x^k\in M_2 (\complex[x])$ for $m\geq 4$  satisfies $P_0 (x\theta(x))=
P_0 (\theta(0)x)$, $P_0 (x^j \theta(x))=0$ for $j\geq 2$, 
$$\sum_{k=0}^{q}(-1)^k S_2 ^{k-q+1}P_k (\theta)=0, \hspace{0.cm} \text{for} \hspace{0.1cm} 0\leq q \leq 1.$$ Then, $[S_2,a_0]=0$, $[S_2,a_1]=[S_2,a_0]$, 
$S_2 a_2 S_2=S_2 a_1$, and $S_2 a_3 S_2=a_2 S_2+ S_2 a_2-a_1$. Writing 
$$a_k =\paren{
	\begin{matrix} 
	r_{k}^{11} & r_{k}^{12} \\
	r_{k}^{21} & r_{k}^{22} 
	\end{matrix}} \mbox{\rm  ,}$$
we have that $r_{0}^{21}=0$, $r_{0}^{22}=r_{0}^{11}$, $r_{1}^{21}=0$, $r_{1}^{22}=r_{1}^{11}=r_{2}^{21}$ and $r_{3}^{21}=r_{2}^{22}+r_{2}^{11}-r_{1}^{12}$. These  equations are exactly those that describe the algebra of the form of Equation~\eqref{motivation} as a sub-algebra of $M_{2}(\complex[x])$.

We show now some properties of the family $\mathcal{P}:= \llave{P_{k}}_{k\in \natu}$.
\begin{lemma}\label{PK0}
	For every $\theta \in M_{N}(\mathbb{K}[x])$,
	\begin{equation}
	P_{k}(\theta)=P_{0}\paren{\frac{\theta^{(k+1)}(0)}{k!}x-\sum_{r=2}^{D}\frac{\theta^{(r+k)}(0)}{(r+k)!}x^r}.
	\end{equation}
\end{lemma}

\begin{demos}
	In fact,
	\begin{equation*}
	P_{0}\paren{\frac{\theta^{(k+1)}(0)}{k!}x-\sum_{r=2}^{D}\frac{\theta^{(r+k)}(0)}{(r+k)!}x^r}
	=\frac{\theta^{(k+1)}(0)}{k!}-\sum_{r=2}^{D}(-1)^r\corch{\frac{\theta^{(r+k)}(0)}{(r+k)!},S^{r-1}}
	\end{equation*}
	\begin{equation*}
		=\frac{\theta^{(k+1)}(0)}{k!}-\sum_{j=k+2}^{k+D}(-1)^{k-j}\corch{\frac{\theta^{(j)}(0)}{j!},S^{j-k-1}}
	=P_{k}(\theta).
	\end{equation*}

\end{demos}
The previous lemma allows us to study the properties of the family $\mathcal{P}=\llave{P_{k}}_{k\in \natu}$ through $P_0$. 

\begin{lemma}[Product Formula for $P_{0}$]\label{Product P0}
	If $\theta_{1}, \theta_{2}\in M_{N}(\mathbb{K}[x])$ then,
	\begin{equation*}
	P_{0}(\theta_{1}\theta_{2})=\sum_{s=0}^{D}\llave{P_{0}(x^s \theta_{1}(x))\frac{\theta_{2}^{(s)}(0)}{s!}+\frac{\theta_{1}^{(s)}(0)}{s!}P_{0}(
		x^s \theta_{2}(x))}-(\theta_{1}\theta_{2})'(0).
	\end{equation*}
\end{lemma}
\begin{demos}
By the definition of $P_0$, 
{\small		\begin{equation*}
	P_{0}(\theta_{1}\theta_{2})= (\theta_{1}\theta_{2})'(0)-\sum_{r=2}^{D}(-1)^r \corch{\frac{(\theta_{1}\theta_{2})^{r}(0)}{r!},S^{r-1}}
	\end{equation*}
	\begin{equation*}
	=\theta_{1}'(0)\theta_{2}(0)+\theta_{1}(0)\theta_{2}'(0)
	-\sum_{r=2}^{D}(-1)^r \corch{\sum_{t=0}^{r}\frac{\theta_{1}^{(t)}(0)}{t!}\frac{\theta_{2}^{(r-t)}(0)}{(r-t)!},S^{r-1}}
    \end{equation*}
    \begin{equation*}
    =\theta_{1}'(0)\theta_{2}(0)+\theta_{1}(0)\theta_{2}'(0)-
    \sum_{r=2}^{D}\sum_{t=0}^{r}(-1)^r \corch{\frac{\theta_{1}^{(t)}(0)}{t!},S^{r-1}}\frac{\theta_{2}^{(r-t)}(0)}{(r-t)!}
       \end{equation*}
    \begin{equation*}
    -\sum_{r=2}^{D}\sum_{t=0}^{r}(-1)^r \frac{\theta_{1}^{(r-t)}(0)}{(r-t)!}\corch{\frac{\theta_{2}^{(t)}(0)}{t!},S^{r-1}}
     =\theta_{1}'(0)\theta_{2}(0)+\theta_{1}(0)\theta_{2}'(0)- 
    \sum_{r=2}^{D}(-1)^r \corch{\frac{\theta_{1}^{(r)}(0)}{r!},S^{r-1}}\theta_{2}(0)
    \end{equation*}
    
    \begin{equation*}
   - \sum_{r=2}^{D}\sum_{t=0}^{r-1}(-1)^r \corch{\frac{\theta_{1}^{(t)}(0)}{t!},S^{r-1}}\frac{\theta_{2}^{(r-t)}(0)}{(r-t)!}
       \end{equation*}
        \begin{equation*}
      -\sum_{r=2}^{D}(-1)^r \theta_{1}(0)\corch{\frac{\theta_{2}^{(r)}(0)}{r!},S^{r-1}}      
      -\sum_{r=2}^{D}\sum_{t=0}^{r-1}(-1)^r \frac{\theta_{1}^{(r-t)}(0)}{(r-t)!}\corch{\frac{\theta_{2}^{(t)}(0)}{t!},S^{r-1}}
    \end{equation*}
     \begin{equation*}
     =\paren{ \theta_{1}'(0)- \sum_{r=2}^{D}(-1)^r \corch{\frac{\theta_{1}^{(r)}(0)}{r!},S^{r-1}} }\theta_{2}(0)
      +\theta_{1}(0)\paren{ \theta_{2}'(0)- \sum_{r=2}^{D}(-1)^r \corch{\frac{\theta_{2}^{(r)}(0)}{r!},S^{r-1}} }
    \end{equation*}
    \begin{equation*}
    - \sum_{r=2}^{D}\sum_{t=0}^{r-1}(-1)^r \corch{\frac{\theta_{1}^{(t)}(0)}{t!},S^{r-1}}\frac{\theta_{2}^{(r-t)}(0)}{(r-t)!}
     -\sum_{r=2}^{D}\sum_{t=0}^{r-1}(-1)^r \frac{\theta_{1}^{(r-t)}(0)}{(r-t)!}\corch{\frac{\theta_{2}^{(t)}(0)}{t!},S^{r-1}}
    \end{equation*}
     \begin{equation*}
     =P_0 (\theta_{1})\theta_{2}(0)+\theta_{1}(0)P_0 (\theta_{2})
      -\sum_{r=2}^{D}(-1)^r \corch{\theta_{1}(0),S^{r-1}} \frac{\theta_{2}^{(r)}(0)}{r!}   -
        \sum_{r=2}^{D}(-1)^r \frac{\theta_{1}^{(r)}(0)}{r!}\corch{\theta_{2}(0),S^{r-1}}
    \end{equation*}
     \begin{equation*}
    - \sum_{r=2}^{D}\sum_{t=1}^{r-1}(-1)^r \corch{\frac{\theta_{1}^{(t)}(0)}{t!},S^{r-1}}\frac{\theta_{2}^{(r-t)}(0)}{(r-t)!}
    -\sum_{r=2}^{D}\sum_{t=1}^{r-1}(-1)^r \frac{\theta_{1}^{(r-t)}(0)}{(r-t)!}\corch{\frac{\theta_{2}^{(t)}(0)}{t!},S^{r-1}}
    \end{equation*}
     \begin{equation*}
    =P_0 (\theta_{1})\theta_{2}(0)+\theta_{1}(0)P_0 (\theta_{2})
    -\sum_{r=2}^{D}(-1)^r \paren{\corch{\theta_{1}(0),S^{r-1}} \frac{\theta_{2}^{(r)}(0)}{r!} + \frac{\theta_{1}^{(r)}(0)}{r!}\corch{\theta_{2}(0),S^{r-1}}}
    \end{equation*}
      \begin{equation*}
    - \sum_{r=2}^{D}\sum_{t=1}^{r-1}(-1)^r \paren{\corch{\frac{\theta_{1}^{(t)}(0)}{t!},S^{r-1}}\frac{\theta_{2}^{(r-t)}(0)}{(r-t)!}+\frac{\theta_{1}^{(r-t)}(0)}{(r-t)!}\corch{\frac{\theta_{2}^{(t)}(0)}{t!},S^{r-1}}}.
    \end{equation*}
   However, 
     \begin{equation*}
     \sum_{r=2}^{D}\sum_{t=1}^{r-1}(-1)^r \paren{\corch{\frac{\theta_{1}^{(t)}(0)}{t!},S^{r-1}}\frac{\theta_{2}^{(r-t)}(0)}{(r-t)!}+\frac{\theta_{1}^{(r-t)}(0)}{(r-t)!}\corch{\frac{\theta_{2}^{(t)}(0)}{t!},S^{r-1}}}
    \end{equation*}
     \begin{equation*}
  =   \sum_{t=1}^{D-1}\sum_{r=t+1}^{D}(-1)^r \paren{\corch{\frac{\theta_{1}^{(t)}(0)}{t!},S^{r-1}}\frac{\theta_{2}^{(r-t)}(0)}{(r-t)!}+\frac{\theta_{1}^{(r-t)}(0)}{(r-t)!}\corch{\frac{\theta_{2}^{(t)}(0)}{t!},S^{r-1}}}
    \end{equation*}
      \begin{equation*}
    =   \sum_{t=1}^{D-1}\sum_{s=1}^{D-t}(-1)^{s+t} \paren{\corch{\frac{\theta_{1}^{(t)}(0)}{t!},S^{s+t-1}}\frac{\theta_{2}^{(s)}(0)}{(s)!}+\frac{\theta_{1}^{(s)}(0)}{(s)!}\corch{\frac{\theta_{2}^{(t)}(0)}{t!},S^{s+t-1}}}
    \end{equation*}
      \begin{equation*}
    =   \sum_{s=1}^{D-1}\sum_{t=1}^{D-s}(-1)^{s+t} \paren{\corch{\frac{\theta_{1}^{(t)}(0)}{t!},S^{s+t-1}}\frac{\theta_{2}^{(s)}(0)}{(s)!}+\frac{\theta_{1}^{(s)}(0)}{(s)!}\corch{\frac{\theta_{2}^{(t)}(0)}{t!},S^{s+t-1}}}
    \end{equation*}
     \begin{equation*}
    =   \sum_{s=1}^{D-1}\sum_{u=s+1}^{D}(-1)^{u} \paren{\corch{\frac{\theta_{1}^{(u-s)}(0)}{(u-s)!},S^{u-1}}\frac{\theta_{2}^{(s)}(0)}{(s)!}+\frac{\theta_{1}^{(s)}(0)}{(s)!}\corch{\frac{\theta_{2}^{(u-s)}(0)}{(u-s)!},S^{u-1}}}
    \end{equation*}
     \begin{equation*}
    = -  \sum_{s=1}^{D-1} \llave{P_{0}\paren{\sum_{u=s+1}^{D}\frac{\theta_{1}^{(u-s)}(0)}{(u-s)!}x^u}\frac{\theta_{2}^{(s)}(0)}{(s)!}+\frac{\theta_{1}^{(s)}(0)}{(s)!}P_0\paren{\sum_{u=s+1}^{D}\frac{\theta_{2}^{(u-s)}(0)}{(u-s)!}x^u}}.
    \end{equation*}
    On the other hand, 
     \begin{equation*}
\sum_{r=2}^{D}(-1)^r \paren{\corch{\theta_{1}(0),S^{r-1}} \frac{\theta_{2}^{(r)}(0)}{r!} + \frac{\theta_{1}^{(r)}(0)}{r!}\corch{\theta_{2}(0),S^{r-1}}}
  \end{equation*}
 \begin{equation*}
=(-1)^D \paren{\corch{\theta_{1}(0),S^{D-1}} \frac{\theta_{2}^{(D)}(0)}{D!} + \frac{\theta_{1}^{(D)}(0)}{D!}\corch{\theta_{2}(0),S^{D-1}}}
    \end{equation*}
     \begin{equation*}
  +  \sum_{r=2}^{D-1}(-1)^r \paren{\corch{\theta_{1}(0),S^{r-1}} \frac{\theta_{2}^{(r)}(0)}{r!} + \frac{\theta_{1}^{(r)}(0)}{r!}\corch{\theta_{2}(0),S^{r-1}}}
    \end{equation*}
    \begin{equation*}
    =(-1)^D \paren{\corch{\theta_{1}(0),S^{D-1}} \frac{\theta_{2}^{(D)}(0)}{r!} + \frac{\theta_{1}^{(D)}(0)}{D!}\corch{\theta_{2}(0),S^{D-1}}}
    \end{equation*}
     \begin{equation*}
    -  \sum_{s=1}^{D-1} \paren{P_0 (\theta_{1}(0)x^s) \frac{\theta_{2}^{(s)}(0)}{s!} + \frac{\theta_{1}^{(s)}(0)}{s!}P_0 (\theta_{2}(0)x^s)}+(\theta_{1}\theta_{2})'(0).
    \end{equation*}
    Then, 
    	\begin{equation*}
    P_{0}(\theta_{1}\theta_{2})=P_0 (\theta_{1})\theta_{2}(0)+\theta_{1}(0)P_0 (\theta_{2})-(\theta_{1}\theta_{2})'(0)
    \end{equation*}
      \begin{equation*}
    +  \sum_{s=1}^{D-1} \llave{P_{0}\paren{\sum_{u=s}^{D}\frac{\theta_{1}^{(u-s)}(0)}{(u-s)!}x^u}\frac{\theta_{2}^{(s)}(0)}{(s)!}+\frac{\theta_{1}^{(s)}(0)}{(s)!}P_0\paren{\sum_{u=s}^{D}\frac{\theta_{2}^{(u-s)}(0)}{(u-s)!}x^u}}
    \end{equation*}
    \begin{equation*}
    -(-1)^D \paren{\corch{\theta_{1}(0),S^{D-1}} \frac{\theta_{2}^{(D)}(0)}{D!} + \frac{\theta_{1}^{(D)}(0)}{D!}\corch{\theta_{2}(0),S^{D-1}}}
    \end{equation*}
    	\begin{equation*}
=\sum_{s=0}^{D}\llave{P_{0}(x^s \theta_{1}(x))\frac{\theta_{2}^{(s)}(0)}{s!}+\frac{\theta_{1}^{(s)}(0)}{s!}P_{0}(
    	x^s \theta_{2}(x))}-(\theta_{1}\theta_{2})'(0).
    \end{equation*}}\qed
\end{demos}

Some remarkable cases of the Lemma \ref{Product P0} are stated in the following corollaries.
\begin{coro}
	If $\theta_{1},\theta_{2}\in M_{N}(\mathbb{K}[x])$ with $\theta_{1}=c\in M_N (\mathbb{K})$ is a constant, then 
	\begin{equation*}
	P_{0}(c\theta_{2})=cP_{0}(\theta_{2})+\sum_{s=2}^{D}P_{0}(cx^s).\frac{\theta_{2}^{(s)}(0)}{s!}.
	\end{equation*}
\end{coro}
\begin{coro}\label{zero1}
If $\theta_{1}(0)=0$, then
\begin{equation*}
P_{0}(\theta_{1}\theta_{2})=P_{0}(\theta_{1})\theta_{2}(0)-\theta_{1}'(0)P_{0}(\theta_{2}(0)x)+\sum_{s=1}^{D}\llave{P_{0}(x^s \theta_{1}(x))\frac{\theta_{2}^{(s)}(0)}{s!}+\frac{\theta_{1}^{(s)}(0)}{s!}P_{0}(
	x^s \theta_{2}(x))}.
\end{equation*} 	
\end{coro}
\begin{coro}\label{zero12}
If $\theta_{1}(0)=\theta_{2}(0)=0$,  then 
\begin{equation*}
P_{0}(\theta_{1}\theta_{2})=\sum_{s=1}^{D}\llave{P_{0}(x^s \theta_{1}(x))\frac{\theta_{2}^{(s)}(0)}{s!}+\frac{\theta_{1}^{(s)}(0)}{s!}P_{0}(
	x^s \theta_{2}(x))}.
\end{equation*}	
\end{coro}
The next lemma tells us that knowledge of any  $P_0$ determines the family $\mathcal{P}=\llave{P_{k}}_{k\in \natu}$.
\begin{lemma}
	For every $\theta \in M_{N}(\mathbb{K}[x])$, we have that
	\begin{equation*}
	P_{0}(\theta)=-\frac{k}{k+1}P_{k}(\theta'(0)x^{k+1})+P_{k}(x^{k}(\theta(x)-\theta(0))).
	\end{equation*}
\end{lemma}
\begin{teor}[Product Formula for $P_{k}$]\label{Product Pk}
	If $\theta_{1}, \theta_{2} \in M_{N}(\mathbb{K}[x])$, then
	\begin{equation*}
	P_{k}(\theta_{1}\theta_{2})=\sum_{t=0}^{k+D}\llave{P_{k}(x^t \theta_{1}(x))\frac{\theta_{2}^{(t)}(0)}{t!}+\frac{\theta_{1}^{(t)}(0)}{t!}P_{k}(
		x^t \theta_{2}(x))}-\frac{(\theta_{1}\theta_{2})^{(k+1)}(0)}{k!}.
	\end{equation*} 
\end{teor}
\begin{demos}
	It is an application of the Lemma \ref{PK0}  and the Lemma \ref{Product P0} .
\end{demos}
\begin{lemma}[Translation]\label{translation}
	For every $\theta\in M_{N}(\mathbb{K}[x])$, $k\geq 0$, we have that 
	\begin{equation}
	P_{k}(\theta) = \left\{
	\begin{array}{ll}
	P_{k}(\frac{\theta^{(k+1-t)}(0)}{(k+1-t)!}x^{k+1})+P_{k-t}(\theta)-\frac{\theta^{(k-t+1)}(0)}{(k-t)!}      & \mbox{\rm if\ }  0\leq t \leq k,\\
	P_{k}(\theta(0)x^{k+1})+P_{0}(x(\theta(x)-\theta(0))) & \mbox{\rm if\ } t=k+1, \\
	P_{0}(x^{t-k}\theta(x))   & \mbox{\rm if\ } t\geq k+2.
	\end{array}
	\right.
	\end{equation}
\end{lemma}
\begin{demos}
	The proof is a straightforward computation.
\end{demos}
We shall now provide the aforementioned description of the Algebra $\mathbb{A} $.
\begin{teor}\label{Gen}
	Let 
\begin{multline*}
 \Gamma :=\Big\{\theta\in M_{N}(\mathbb{K}[x])\mid P_{0}(x\theta(x))=P_{0}(\theta(0)x), P_{0}(x^{j}\theta(x))=0, j\geq 2 ,
  \\
 \sum_{k=0}^{q}(-1)^k S^{k+D-q-1} P_{k}(\theta)=0, 0\leq q \leq D-1 \Big\}. 
 \end{multline*}
	Then, $\Gamma=\mathbb{A}$.
	
	Moreover,
	for each $\theta$ we have an explicit expression for the operator $\mathcal{B}$. 
\end{teor}
Before proving the theorem, we study the relations defining the algebra $\Gamma$. They are given in the following result: 
\begin{prop}\label{rel}
	The algebra $\Gamma$ is the subset of $\theta \in M_{N}(\mathbb{K}[x])$ such that 
	
	\begin{equation}\label{{1fam}}
	\sum_{j=0}^{q}(-1)^{q-j-D}\corch{S^{D-q+j-1}, \frac{\theta^{(j)}(0)}{j!}}=0,
	\end{equation}
	\begin{equation}\label{2fam}
	\sum_{j=0}^{q}(-1)^{q-j-D+1}S^{j+D-q-1}P_{j}(\theta)=0,
	\end{equation}
	for $0\leq q \leq D-1$.
\end{prop}
\begin{demos}
We notice that $\Gamma$ is defined by two relations:
\begin{equation}\label{eq1}
  P_{0}(x\theta(x))=P_{0}(\theta(0)x), P_{0}(x^{j}\theta(x))=0, \; j\geq 2    
\end{equation}
and Equation~\eqref{2fam} (after a trivial change of the summation variable).
	Equation~\eqref{eq1} is equivalent to 
		\begin{equation*}
	\sum_{r=D-q}^{D}(-1)^{r}\corch{S^{r-1},\frac{\theta^{(r-D+q)}(0)}{(r-D+q)!}}=0,
	\end{equation*}
	for $0\leq q \leq D-1$.
	If $q=D-1$, then 
		\begin{equation*}
0=	\sum_{r=2}^{D}(-1)^{r}\corch{S^{r-1},\frac{\theta^{(r-1)}(0)}{(r-1)!}}=P_{0}\paren{\sum_{r=2}^{D}\frac{\theta^{(r-1)}(0)}{(r-1)!}x^r}
	\end{equation*}
	\begin{equation*}
	=P_{0}\paren{\sum_{r=1}^{D-1}\frac{\theta^{(r)}(0)}{r!}x^{r+1}}=P_{0}(x(\theta(x)-\theta(0))).
	\end{equation*}
	In other words, $P_{0}(x\theta(x))=P_{0}(x\theta(0))$.
	
	If $0\leq q \leq D-2$, then $2 \leq D-q$ and 
	\begin{equation*}
	0= \sum_{r=D-q}^{D}(-1)^{r}\corch{S^{r-1},\frac{\theta^{(r-D+q)}(0)}{(r-D+q)!}}=P_{0}\paren{\sum_{r=D-q}^{D}\frac{\theta^{(r-D-q)}(0)}{(r-D-q)!}x^r}
	\end{equation*}
	\begin{equation*}
	=P_{0}\paren{\sum_{j=0}^{q}\frac{\theta^{(j)}(0)}{j!}x^{j+D-q}}=P_{0}\paren{x^{D-q}\paren{\sum_{j=0}^{q}\frac{\theta^{(j)}(0)}{j!}x^{j}}}
	=P_{0}(x^{D-q}\theta(x))
	\end{equation*}
	for $0\leq q \leq D-2$. Thus,	$P_{0}(x^j\theta(x))=0$ for $j\geq 2$.
\end{demos}
Now let us prove the theorem.
\begin{demos}
	We shall break the proof in different steps.\\
	Step 1: The set $\Gamma$ is an algebra.
	
	Clearly,\hspace{0.1cm} $\Gamma$ is a vector space since $P_{k}$ is linear for all $0\leq k \leq D-1$. 
	
	If $\theta_{1}, \theta_{2}\in \Gamma$, then 
	 \begin{equation*}
	 P_{0}(x\theta_{i}(x))=P_{0}(\theta_{i}(0)x),\; 
	 P_{0}(x^{j}\theta_{i}(x))=0, \; 
	 \sum_{k=0}^{q}(-1)^k S^{k+D-q-1} P_{k}(\theta_{i})=0,
	\end{equation*}
	for $j\geq 2,\hspace{0.1cm} 0\leq q \leq D-1,\hspace{0.1cm} i=1,2$.
	
	Note that, using Corollary \ref{zero1} and $P_0(x\theta_{1}(x))(0)=0$, we obtain
	\begin{equation*}
	P_{0}(x\theta_{1}(x)\theta_{2}(x))=P_{0}((x\theta_{1}(x))\theta_{2}(x))
	=P_{0}(\theta_{1}(0)x)\theta_{2}(0)-\theta_{1}(0)P_{0}(\theta_{2}(0)x)+
	(x\theta_{1})'(0)P_{0}(x\theta_{2}(x))
    \end{equation*}
    \begin{equation*}	
	=P_{0}(\theta_{1}(0)x)\theta_{2}(0)
	-\theta_{1}(0)P_{0}(\theta_{2}(0)x)+\theta_{1}(0)P_{0}(\theta_{2}(0)x)
	=\theta_{1}(0)\theta_{2}(0)= P_{0}(\theta_{1}(0)\theta_{2}(0)x).
	\end{equation*}
	If $j \geq 2$, then $P_0(x^{j-1}\theta_{1}(x))(0)=0, P_0(x\theta_{2}(x))(0)=0$. Using  Corollary \ref{zero12} we obtain 
		\begin{equation*}
	P_{0}(x^{j}\theta_{1}(x)\theta_{2}(x))=P_{0}((x^{j-1}\theta_{1}(x))(x\theta_{2}(x)))
	\end{equation*}
	\begin{equation*}
	=\sum_{s=1}^{D}\llave{P_{0}(x^{j+s-1} \theta_{1}(x))\frac{(x\theta_{2})^{(s)}(0)}{s!}+\frac{(x^{j-1}\theta_{1})^{(s)}(0)}{s!}P_{0}(
		x^{s+1} \theta_{2}(x))}=0
	\end{equation*}
	Lemma~\ref{translation} (Translation)  implies 
	 \begin{equation*}
	 P_{k}(x^{k+1}\theta_{i}(x))=P_{k}(\theta_{i}(0)x^{k+1}),
	 P_{k}(x^t\theta_{i}(x))=P_{0}(x^{t-k}\theta_{i}(x))=0, 
	 \end{equation*}
	 for $t\geq k+2$,\hspace{0.1cm} $i=1,2$.
	 
	 Using  Theorem \ref{Product Pk}  (Product Formula for $P_{k}$) we have:
	 
	 \begin{equation*}
	 P_{k}(\theta_{1}\theta_{2})=\sum_{t=0}^{k+D}\llave{P_{k}(x^t \theta_{1}(x))\frac{\theta_{2}^{(t)}(0)}{t!}+\frac{\theta_{1}^{(t)}(0)}{t!}P_{k}(
	 	x^t \theta_{2}(x))}-\frac{(\theta_{1}\theta_{2})^{(k+1)}(0)}{k!}
	 \end{equation*}
	 \begin{equation*}
	 =\sum_{t=0}^{k}\llave{P_{k}(x^t \theta_{1}(x))\frac{\theta_{2}^{(t)}(0)}{t!}+\frac{\theta_{1}^{(t)}(0)}{t!}P_{k}(
	 	x^t \theta_{2}(x))}+P_{k}( \theta_{1}(0)x^{k+1})\frac{\theta_{2}^{(k+1)}(0)}{(k+1)!}
 	\end{equation*}
 	\begin{equation*}
 	+\frac{\theta_{1}^{(k+1)}(0)}{(k+1)!}P_{k}(
 	\theta_{2}(0)x^{k+1})-\frac{(\theta_{1}\theta_{2})^{(k+1)}(0)}{k!}
	 \end{equation*}
	 \begin{equation*}
	 =P_{k}(\theta_{1})\theta_{2}(0)+\theta_{2}(0)P_{k}(\theta_{2})
	  \end{equation*}
	 \begin{equation*}
	 + \sum_{t=0}^{k}\llave{P_{k}(x^t \theta_{1}(x))\frac{\theta_{2}^{(t)}(0)}{t!}+\frac{\theta_{1}^{(t)}(0)}{t!}P_{k}(
	 	x^t \theta_{2}(x))-(k+1)\frac{\theta_{1}^{(t)}(0)}{t!}\frac{\theta_{2}^{(k+1-t)}(0)}{(k+1-t)!}}.
	 \end{equation*}
	 Thus, for $0\leq q \leq D-1$ we have 
	 \begin{equation*}
	 \sum_{k=0}^{q}(-1)^k S^{k+D-q-1}P_{k}(\theta_{1}\theta_{2})
	 = \sum_{k=0}^{q}(-1)^k S^{k+D-q-1}\Bigg[ P_{k}(\theta_{1})\theta_{2}(0)+\theta_{2}(0)P_{k}(\theta_{2}) 
	 \end{equation*}
	 \begin{equation*}
	 +	 	\left.\sum_{t=0}^{k}\llave{P_{k}(x^t \theta_{1}(x))\frac{\theta_{2}^{(t)}(0)}{t!}+\frac{\theta_{1}^{(t)}(0)}{t!}P_{k}(
	 		x^t \theta_{2}(x))-(k+1)\frac{\theta_{1}^{(t)}(0)}{t!}\frac{\theta_{2}^{(k+1-t)}(0)}{(k+1-t)!}}\right]
	 \end{equation*}
	\begin{equation*}
	=\sum_{k=0}^{q}(-1)^{k} \corch{S^{k+D-q-1},\theta_{1}(0)}P_{k}(\theta_{2})
	\end{equation*}
	\begin{equation*}
	+
	\sum_{k=0}^{q}\sum_{t=1}^{k}(-1)^{k} S^{k+D-q-1}\llave{(k+1)\frac{\theta_{1}^{(k+1-t)}(0)}{(k+1-t)!}+P_{k-t}(\theta_{1})-
		\frac{\theta_{1}^{(k+1-t)}(0)}{(k-t)!}}\frac{\theta_{2}^{(t)}(0)}{t!}
	\end{equation*}
	\begin{equation*}
	+\sum_{k=0}^{q}\sum_{t=1}^{k}(-1)^{k} S^{k+D-q-1}\frac{\theta_{1}^{(t)}(0)}{t!}
	\llave{(k+1)\frac{\theta_{2}^{(k+1-t)}(0)}{(k+1-t)!}+P_{k-t}(\theta_{2})-
		\frac{\theta_{2}^{(k+1-t)}(0)}{(k-t)!}}
	\end{equation*}
	\begin{equation*}
	-\sum_{k=0}^{q}\sum_{t=1}^{k}(k+1)(-1)^{k} S^{k+D-q-1}\frac{\theta_{1}^{(t)}(0)}{t!}\frac{\theta_{2}^{(k+1-t)}(0)}{(k-t+1)!}
	\end{equation*}
	\begin{equation*}
		=\sum_{k=0}^{q}(-1)^{k} \corch{S^{k+D-q-1},\theta_{1}(0)}P_{k}(\theta_{2})
	\end{equation*}
	\begin{equation*}
	+
	\sum_{k=0}^{q}\sum_{t=1}^{k}(-1)^{k} S^{k+D-q-1}\llave{P_{k-t}(\theta_{1})-
		\frac{\theta_{1}^{(k+1-t)}(0)}{(k-t)!}}\frac{\theta_{2}^{(t)}(0)}{t!}
	\end{equation*}
	\begin{equation*}
	+\sum_{k=0}^{q}\sum_{t=1}^{k}(-1)^{k} S^{k+D-q-1}\frac{\theta_{1}^{(t)}(0)}{t!}
	\llave{P_{k-t}(\theta_{2})-
		\frac{\theta_{2}^{(k+1-t)}(0)}{(k-t)!}}
	\end{equation*}
	\begin{equation*}
	+\sum_{k=0}^{q}\sum_{t=1}^{k}(k+1)(-1)^{k} S^{k+D-q-1}\frac{\theta_{1}^{(t)}(0)}{t!}\frac{\theta_{2}^{(k+1-t)}(0)}{(k-t+1)!}.
	\end{equation*}
	
However,
	\begin{equation*}
	\sum_{k=0}^{q}\sum_{t=1}^{k}(-1)^{k} S^{k+D-q-1}\frac{\theta_{1}^{(t)}(0)}{t!}
	\llave{P_{k-t}(\theta_{2})-
		\frac{\theta_{2}^{(k+1-t)}(0)}{(k-t)!}}
	\end{equation*}
	\begin{equation*}
	=\sum_{s=0}^{q-1}\sum_{k=s+1}^{q}(-1)^{k} S^{k+D-q-1}\frac{\theta_{1}^{(k-s)}(0)}{(k-s)!}\llave{P_{s}(\theta_{2})-
		\frac{\theta_{2}^{(s+1)}(0)}{s!}}
	\end{equation*}
	\begin{equation*}
		=\sum_{s=0}^{q-1}(-1)^s\paren{\sum_{j=1}^{q-s}(-1)^{j} S^{j+s+D-q-1}\frac{\theta_{1}^{j}(0)}{j!}}\llave{P_{s}(\theta_{2})-
		\frac{\theta_{2}^{(s+1)}(0)}{s!}}
	\end{equation*}
	\begin{equation*}
	=\sum_{s=0}^{q-1}(-1)^s\paren{\sum_{j=1}^{q-s}(-1)^{j} \frac{\theta_{1}^{j}(0)}{j!}S^{j+s+D-q-1}-\corch{S^{D-q+s-1},\theta_{1}(0)}}\llave{P_{s}(\theta_{2})-
		\frac{\theta_{2}^{(s+1)}(0)}{s!}}.
	\end{equation*}
	After a few simple calculations this term is equal to:
		\begin{equation*}
	-\sum_{s=0}^{q-1}(-1)^s \corch{S^{D-q+s-1},\theta_{1}(0)}P_{s}(\theta_{2}) 
	-\sum_{k=0}^{q}\sum_{s=0}^{k-1}(-1)^{k} S^{k+D-q-1}\frac{\theta_{1}^{(s+1)}(0)}{(s+1)!}\frac{\theta_{2}^{(k-s)}(0)}{(k-s-1)!}.
	\end{equation*}
	Similarly, we can see that 
	\begin{equation*}
	\sum_{k=0}^{q}\sum_{t=1}^{k}(-1)^{k} S^{k+D-q-1}\llave{P_{k-t}(\theta_{1})-
		\frac{\theta_{1}^{(k+1-t)}(0)}{(k-t)!}}\frac{\theta_{2}^{(t)}(0)}{t!}
	\end{equation*}
	\begin{equation*}
	=-\sum_{k=0}^{q}\sum_{s=0}^{k-1}(-1)^{k} S^{k+D-q-1}\frac{\theta_{1}^{(s+1)}(0)}{s!}\frac{\theta_{2}^{(k-s)}(0)}{(k-s)!}.
	\end{equation*}
	 Thus,
	\begin{equation*}
	\sum_{k=0}^{q}(-1)^k S^{k+D-q-1}P_{k}(\theta_{1}\theta_{2})
	=  \sum_{s=0}^{q}(-1)^s \corch{S^{D-q+s-1},\theta_{1}(0)}P_{s}(\theta_{2}) 
	\end{equation*}
		\begin{equation*}
		-\sum_{k=0}^{q}\sum_{s=0}^{k-1}(-1)^{k} S^{k+D-q-1}\frac{\theta_{1}^{(s+1)}(0)}{s!}\frac{\theta_{2}^{(k-s)}(0)}{(k-s)!}
		\end{equation*}
		
		\begin{equation*}
	-\sum_{s=0}^{q-1}(-1)^s \corch{S^{D-q+s-1},\theta_{1}(0)}P_{s}(\theta_{2}) 
	-\sum_{k=0}^{q}\sum_{s=0}^{k-1}(-1)^{k} S^{k+D-q-1}\frac{\theta_{1}^{(s+1)}(0)}{(s+1)!}\frac{\theta_{2}^{(k-s)}(0)}{(k-s-1)!}
	\end{equation*}
	 \begin{equation*}
	+\sum_{k=0}^{q}\sum_{s=0}^{k-1}(-1)^{k}(k+1) S^{k+D-q-1}\frac{\theta_{1}^{(s+1)}(0)}{(s+1)!}\frac{\theta_{2}^{(k-s)}(0)}{(k-s)!}
	\end{equation*}
	\begin{equation*}
	=(-1)^q \corch{S^{D-1},\theta_{1}(0)}P_{k}(\theta_{2})
	-\sum_{k=0}^{q}\sum_{s=0}^{k-1}(-1)^{k} (s+1)S^{k+D-q-1}\frac{\theta_{1}^{(s+1)}(0)}{s!}\frac{\theta_{2}^{(k-s)}(0)}{(k-s)!}
	\end{equation*}
	\begin{equation*}
-\sum_{k=0}^{q}\sum_{s=0}^{k-1}(-1)^{k}(k-s) S^{k+D-q-1}\frac{\theta_{1}^{(s+1)}(0)}{(s+1)!}\frac{\theta_{2}^{(k-s)}(0)}{(k-s)!}
	\end{equation*}
	\begin{equation*}
	+\sum_{k=0}^{q}\sum_{s=0}^{k-1}(-1)^{k}(k+1) S^{k+D-q-1}\frac{\theta_{1}^{(s+1)}(0)}{(s+1)!}\frac{\theta_{2}^{(k-s)}(0)}{(k-s)!}
	\end{equation*}
	\begin{equation*}
	=(-1)^{q}\corch{S^{D-1},\theta_{1}(0)}P_{q}(\theta_{2}).
	\end{equation*}
	Nevertheless, $\corch{S^{D-1},\theta_{1}(0)}=0$ and we obtain $\sum_{k=0}^{q}(-1)^k S^{k+D-q-1}P_{k}(\theta_{1}\theta_{2})=0$ for $0\leq q \leq D-1$. Therefore, $\theta_{1}\theta_{2}\in \Gamma$. Thus, $\Gamma$ is an algebra.
	
	Step 2:  Since $\Gamma$ contains the ideal $x^{2D}M_{N}(\mathbb{K}[x])$, if we define $E=\oplus_{j=0}^{2D-1}M_{N} (\mathbb{K}[x])_{j}\cap \Gamma$, then we have this step.

    Step 3: The algebra $\Gamma$ is generated by $E$, i.e., $\mathbb{K}\cdot \langle E\rangle=\Gamma$.	

    This step follows  applying  Theorem \ref{example} and since $E$ contains the elements mentioned in that theorem. 	
    
    Step 4:  The inclusion  $\mathbb{A}\cap \oplus_{j=0}^{2D-1} M_{N}(\mathbb{K}[x])_{j}\subset E$.
    
    Let $\theta\in \mathbb{A}\cap \oplus_{j=0}^{2D-1} M_{N}(\mathbb{K}[x])_{j}$ then there exists $\mathcal{B}=\mathcal{B}(z,\partial_{z})$ such that  $(\psi \mathcal{B})(x,z)=\theta(x)\psi(x,z)$. We write  $\theta(x)=\sum_{j=0}^{2D-1}a_{j}x^{j}$. After a few simple computations we obtain that: 
    \begin{equation}\label{bispectral pair 1}
    \mathcal{B}=\sum_{j=0}^{2D-1}\partial_{z}^{j}\cdot\paren{a_{k}+\sum_{l=1}^{2D-1-j}\frac{(-1)^l}{z^l}\sum_{r=j+l-1}^{2D-1}(\mu^{l-1})_{jr}S^{r-j-l+1}P_{r}(\theta)}.
    \end{equation}
   With $\mu \in M_{2D}(\mathbb{K}[x])$ given by 
   \begin{equation*}
    \mu_{rj}  = \left\{
   \begin{array}{ll}
      (-1)^{r-j}   & \mbox{\rm if\ } r+2\leq j \leq \min\llave{r+D,2D-1}, \\
   r & \mbox{\rm if\ } j=r+1, \\
   0   & \mbox{\rm  if  otherwise.}  
   \end{array}
   \right.
   \end{equation*}
  Furthermore,
  \begin{equation*}
  e^{-xz}(\psi \mathcal{B} -\theta \psi)=
  x^{-D} \left(\sum_{q=0}^{D-1}\left\{\sum_{j=0}^{q}(-1)^{q-j-D}\corch{S^{D-q+j-1},a_{j}} \right.\right.
  		\end{equation*}
  		\begin{equation*}
  		\left. \left.+\sum_{j=0}^{q}\sum_{l=1}^{2D-1-j}\frac{(-1)^{q-j-D+l}}{z^{l}}
  		\sum_{r=j+l-1}^{2D-1}(\mu^{l-1})_{jr}S^{r+D-q-l}P_{r}(\theta)\right\}\right)x^q
  \end{equation*}
  \begin{equation*}
  = x^{-D} \left(\sum_{q=0}^{D-1}\left\{\sum_{j=0}^{q}(-1)^{q-j-D}\corch{S^{D-q+j-1},a_{j}} \right.\right.
  \end{equation*}
  \begin{equation*}
  \left. \left.+\sum_{l=1}^{2D-2-q}\paren{\sum_{j=0}^{q}(-1)^{q-j-D+l-1}
  \sum_{r=j+l-1}^{2D-1}(\mu^{l-1})_{jr}S^{r+D-q-l}P_{r}(\theta)}\frac{1}{z^{l}}\right. \right.
  \end{equation*}
   \begin{equation*}
  \left. \left.+\sum_{l=2D-1-q}^{2D-1}\paren{\sum_{j=0}^{2D-l-1 }(-1)^{q-j-D+l}
  \sum_{r=j+l-1}^{2D-1}(\mu^{l-1})_{jr}S^{r+D-q-l}P_{r}(\theta)}\frac{1}{z^{l}}\right\} x^q\right).
  \end{equation*}
  However, 
\begin{equation*}
\sum_{j=0}^{q}(-1)^{q-j-D+l}
\sum_{r=j+l-1}^{2D-1}(\mu^{l-1})_{jr}S^{r+D-q-l}P_{r}(\theta)
\end{equation*}
\begin{equation*}
=\sum_{r=l-1}^{q+l-1}(-1)^{q-D+l}\paren{
\sum_{j=0}^{r-l+1}(-1)^j (\mu^{l-1})_{jr}}S^{r+D-q-l}P_{r}(\theta),
\end{equation*}
for $ 1\leq l \leq 2D-2-q$,\hspace{0.1cm} and 
\begin{equation*}
\sum_{j=0}^{2D-l-1 }(-1)^{q-j-D+l}
\sum_{r=j+l-1}^{2D-1}(\mu^{l-1})_{jr}S^{r+D-q-l}P_{r}(\theta)
\end{equation*}
\begin{equation*}
=\sum_{r=l-1}^{2D-2}(-1)^{q-D+l}\paren{
	\sum_{j=0}^{r-l+1}(-1)^j (\mu^{l-1})_{jr}}S^{r+D-q-l}P_{r}(\theta),
\end{equation*}
for $2D-1-q \leq l \leq 2D-1$.

Note that $r$-th component of $v\mu$ is given by $(v\mu)_{r}=\sum_{j=0}^{2D-1}(-1)^{j}\mu_{jr}=\sum_{j=0}^{2D-1}(-1)^{j}\mu_{jr}
=(-1)^{r-1}(r-1)+\sum_{j=0}^{r-2}(-1)^j (-1)^{r-j}=0$ for $v\in \mathbb{K}^{2D}$ defined by  $v_{j}=(-1)^{j}$ and $0\leq j,r \leq 2D-1$. Therefore,  $v\mu=0$.
Clearly this implies $\sum_{j=0}^{r-l+1}(-1)^j (\mu^{l-1})_{jr}= (v\mu^{l-1})_{r}=0$ for $l \geq 2$.

Therefore, 

 \begin{equation*}
e^{-xz}(\psi \mathcal{B} -\theta \psi)=
x^{-D} \left(\sum_{q=0}^{D-1}\left\{\sum_{j=0}^{q}(-1)^{q-j-D}\corch{S^{D-q+j-1},a_{j}} \right.\right.
\end{equation*}
\begin{equation*}
\left. \left.+\frac{1}{z}\sum_{j=0}^{q}(-1)^{q-j-D+1}S^{j+D-q-1}P_{j}(\theta)\right\}x^{q}\right).
\end{equation*}

Since $\theta\in \mathbb{A}$ we have 
\begin{equation*}
\sum_{j=0}^{q}(-1)^{q-j-D}\corch{S^{D-q+j-1},\frac{\theta^{(j)}(0)}{j!}}=0,
\end{equation*}
 \begin{equation*}
 \sum_{j=0}^{q}(-1)^{q-j-D+1}S^{j+D-q-1}P_{j}(\theta)=0,
 \end{equation*}
 for $0\leq q \leq D-1$. Thus,\hspace{0.1cm} $\theta\in E$.
    
    Step 5: The inclusion $E\subset \mathbb{A}\cap \oplus_{j=0}^{2D-1} M_{N}(\mathbb{K}[x])_{j}$. 
    
    By the previous step we have Equation~\eqref{bispectral pair 1} valid for every $\theta\in E$ and using  Proposition \ref{rel} we obtain that $(\psi \mathcal{B}) (x,z)=\theta(x) \psi(x,z)$. Then, $\theta \in 
    \mathbb{A}\cap \oplus_{j=0}^{2D-1} M_{N}(\mathbb{K}[x])_{j}$.
    
    Furthermore, we have an explicit expression for the operator $\mathcal{B}$. 
    
    If $\theta(x)=\sum_{j=0}^{M}a_{j}x^{j}\in \Gamma$,  then 
    \begin{equation}
    \mathcal{B}=\sum_{j=0}^{M}\partial_{z}^{j}\cdot\paren{a_{k}+\sum_{l=1}^{M-j}\frac{(-1)^l}{z^l}\sum_{r=j+l-1}^{M}(\mu^{l-1})_{jr}S^{r-j-l+1}P_{r}(\theta)}
    \end{equation}
    with $\mu \in M_{M+1}(\mathbb{K}[x])$ given by 
   \begin{equation}
   \mu_{rj}  = \left\{
   \begin{array}{ll}
   (-1)^{r-j}   & \mbox{\rm if\ } r+2\leq j \leq \min\llave{r+D,M}, \\
   r & \mbox{\rm if\ } j=r+1, \\
   0   & \mbox{\rm if\ otherwise. } 
   \end{array}
   \right.
   \end{equation}
    satisfies $(\psi \mathcal{B}) (x,z)=\theta(x) \psi(x,z)$.
    $\qed$

\end{demos}

In particular, Theorem~\ref{Gen} implies that the $\mathbb{K}$-algebra $\Gamma$ is not trivial.

A remarkable property of this family of algebras is the existence of a Pierce decomposition whose definition we shall now recall. 
\begin{defin}
Let $R$ be a noncommutative ring with unit. We say that a set of elements $r_1, ... ,r_n \in R$ is  Pierce decomposition of $R$ if $1=\sum_{j=1}^{n}r_j$ and $r_i r_j =\delta_{ij}$ for all $1\leq i,j \leq n$.
\end{defin}

See \cite{Pidecom} for more information on the Pierce decomposition. The next definition presents a Pierce decomposition of the algebra $\mathbb{A}$. 

\begin{defin}
	Define $\alpha_{k}(x)=e_{kk}+\sum_{j=1}^{N-1}a_{kj}x^j \in M_{N}(\mathbb{K}[x])$ for $2\leq k \leq N-1$, $N\geq 3$ with 
	\begin{equation*}
	a_{kj}=(-1)^{j+1}(\delta_{k,j+1}e_{k1}+\delta_{k,N-j}e_{Nk})+(-1)^j \delta_{N,j+1}e_{N1}
	=(-1)^{j+1}(\delta_{j,k-1}e_{k1}+\delta_{j,N-k}e_{Nk})+(-1)^j \delta_{j,N-1}e_{N1}
	\end{equation*}
	for $1\leq j \leq N-1$.
\end{defin}

In the previous definition we have two cases:
\begin{itemize}
	\item If $N$ is even and the numbers $k-1, N-k, N-1$ are different, then 
	 \begin{equation}
	a_{kj}  = \left\{
	\begin{array}{ll}
	(-1)^{k} e_{k1}   & \mbox{\rm if\ } j=k-1, \\
	(-1)^{N-k+1}e_{Nk} & \mbox{\rm if\ } j=N-k, \\
		(-1)^{N-1}e_{N1} & \mbox{\rm if\ } j=N-1, \\
		e_{kk} & \mbox{\rm if\ } j=0, \\
	0   & \mbox{\rm otherwise. } 
	\end{array}
	\right.
	\end{equation}
	\item If $N$ is odd
	\begin{itemize}
			\item If $k\neq \frac{N+1}{2}$, then $k-1, N-k, N-1$  are different:
		\begin{equation}
		a_{kj}  = \left\{
		\begin{array}{ll}
		(-1)^{k} e_{k1}   & \mbox{\rm if\ } j=k-1 ,\\
		(-1)^{N-k+1}e_{Nk} & \mbox{\rm if\ } j=N-k, \\
		(-1)^{N-1}e_{N1} & \mbox{\rm if\ } j=N-1, \\
		e_{kk} & \mbox{\rm if\ } j=0, \\
		0   & \mbox{\rm otherwise. } 
		\end{array}
		\right.
		\end{equation}
		
		\item  If $k=\frac{N+1}{2}$, then
		\begin{equation}
		a_{\frac{N+1}{2}j}  = \left\{
		\begin{array}{ll}
		(-1)^{\frac{N+1}{2}}\paren{e_{\frac{N+1}{2},1}+e_{N,\frac{N+1}{2}}}   & \mbox{\rm if\ } j=\frac{N-1}{2}, \\
		(-1)^{N-1}e_{N1} & \mbox{\rm if\ } j=N-1, \\
		e_{kk} & \mbox{\rm if\ } j=0, \\
		0   & \mbox{\rm otherwise. } 
		\end{array}
		\right.
		\end{equation}
	
	\end{itemize}
	
\end{itemize}
The following lemma relates these elements with the family $\llave{P_k}_{k\in \natu}$. Its importance is that it shall be used to prove that the elements $\alpha_k 's$ satisfy the second family of relations that defines $\mathbb{A}$.

\begin{lemma}\label{pierce}
	\begin{itemize}
		\item 	If $k<\frac{N+1}{2}$, then $k-1<N-k<N-1$ and 
		\begin{equation}
	P_l (\alpha_k) = \left\{
	\begin{array}{ll}
	(-1)^{l} \paren{e_{k,k-l-1}-e_{k+l+1,k}}   & \mbox{\rm if\ } 0\leq l \leq k-3, \\
	(-1)^{k}ke_{k1}+(-1)^{k+1}e_{2k-1,k} & \mbox{\rm if\ } l=k-2, \\
	(-1)^{l}\paren{e_{l+2,1}-e_{k+l+1,k}} & \mbox{\rm if\ } k-1\leq l \leq N-k-2,  \\
	(N-k-1)(-1)^{N-k-1}e_{Nk}+(-1)^{N-k+1}e_{N-k+1,1} & \mbox{\rm if\ } l=N-k-1, \\
	(-1)^{l+1}\paren{e_{N,N-l-1}-e_{l+2,1}}  & \mbox{\rm if\ }  N-k\leq l \leq N-3, \\
	(N-1)(-1)^{N-1}e_{N1} & \mbox{\rm if\ }  N-k\leq l \leq N-3. \\
	\end{array}
	\right.
	\end{equation}
	\item 	If $k=\frac{N+1}{2}$, then $N-k=k-1=\frac{N-1}{2}<N-1$ and 
	\begin{equation}
	P_l (\alpha_k) = \left\{
	\begin{array}{ll}
	(-1)^{l} \paren{e_{\frac{N+1}{2},\frac{N-1}{2}-l}-e_{\frac{N+3}{2}+l,\frac{N+1}{2}}}   & \mbox{\rm if\ } 0\leq l \leq k-3=N-k-2=\frac{N-5}{2}, \\
	\frac{N-3}{2}(-1)^{\frac{N+1}{2}}e_{N,\frac{N+1}{2}}+(-1)^{\frac{N+1}{2}}e_{\frac{N+1}{2},1} & \mbox{\rm if\ } l=k-2=N-k-1=\frac{N-3}{2},\\
	(-1)^{l+1}\paren{e_{N,N-l-1}-e_{l+2,1}}  & \mbox{\rm if\ }  N-k\leq l \leq N-3, \\
	(N-1)(-1)^{N-1}e_{N1} & \mbox{\rm if\ }   l = N-2. \\
	\end{array}
	\right.
	\end{equation}
	\item 	If $k>\frac{N+1}{2}$, then $N-k<k-1<N-1$ and 
	\begin{equation}
	P_l (\alpha_k) = \left\{
	\begin{array}{ll}
	(-1)^{l} \paren{e_{k,k-l-1}-e_{k+l+1,k}}   & \mbox{\rm if\ } 0\leq l \leq N-k-2,\\
	(N-k-1)(-1)^{N-k+1}e_{Nk}+(-1)^{N-k+1}e_{k,2k-N} & \mbox{\rm if\ } l=N-k-1, \\
	(-1)^{l}\paren{e_{k,k-l-1}-e_{N,N-l-1}} & \mbox{\rm if\ } N-k\leq l \leq k-3,  \\
	k(-1)^{k}e_{k1}+(-1)^{k+1}e_{N,N-k+1} & \mbox{\rm if\ } l=k-2, \\
	(-1)^{l+1}\paren{e_{N,N-l-1}-e_{l+2,1}}  & \mbox{\rm if\ }  k-1\leq l \leq N-3, \\
	(N-1)(-1)^{N-1}e_{N1} & \mbox{\rm if\ }   l = N-2. \\
	\end{array}
	\right.
	\end{equation}
	\end{itemize}

\end{lemma}
\begin{demos}
	The proof is a straightforward. \qed
\end{demos}

Now we prove that this family is contained in $\mathbb{A}$.
\begin{teor}\label{pierce A}
	For $N\geq 3$ we have $\llave{\alpha_k}_{1\leq k \leq N-1}\subset \mathbb{A}$.
\end{teor}
Before proving this theorem we have a handy remark.
\begin{remark}
	We shall adopt the convenient convention that $\sum_{i\in \emptyset}x_i =0$. Define $e_{ij}=0$ if $i$ or $j$ is outside the set $\llave{1, ... ,N}$.
\end{remark}
\begin{demos}
	We first verify the first family of relations. Pick  $2\leq k \leq N-1$, then 
	\begin{equation*}
	P_{0}(x\alpha_{k}(x))=e_{kk}-(-1)^{k}\corch{(-1)^{k}e_{k1},S_{N}^{k-1}}
	-(-1)^{N-k+1}\corch{(-1)^{N-k+1}e_{Nk},S_{N}^{N-k}}
	\end{equation*}
	\begin{equation*}
		-(-1)^{N}\corch{(-1)^{N-1}e_{N1},S_{N}^{N-1}}
	=e_{kk}-\corch{e_{k1},S_{N}^{k-1}}-\corch{e_{Nk},S_{N}^{N-k}}-\corch{e_{N1},S_{N}^{N-1}}
	=e_{kk}-(e_{kk}-e_{11})
		\end{equation*}
	\begin{equation*}
-(e_{NN}-e_{kk})+(e_{NN}-e_{11})	=e_{kk}=P_{0}(\alpha_{k}(0)x).
	\end{equation*}
		If $N$ is odd 
	\begin{equation*}
	P_{0}\paren{x\alpha_{\frac{N+1}{2}}(x)}=e_{\frac{N+1}{2},\frac{N+1}{2}}-(-1)^{\frac{N+1}{2}}\corch{(-1)^{\frac{N+1}{2}}\paren{e_{\frac{N+1}{2},1}+e_{N,\frac{N+1}{2}}},S_{N}^{\frac{N-1}{2}}}
	-(-1)^{N}\corch{(-1)^{N-1}e_{N1},S_{N}^{N-1}}
	\end{equation*}
	\begin{equation*}
	=e_{\frac{N+1}{2},\frac{N+1}{2}}-\paren{e_{\frac{N+1}{2},\frac{N+1}{2}}+e_{NN}-\paren{e_{11}+e_{\frac{N+1}{2}},\frac{N+1}{2}}}+e_{NN}-e_{11}
	=e_{\frac{N+1}{2},\frac{N+1}{2}}=P_{0}\paren{\alpha_{\frac{N+1}{2}}(0)x}.
	\end{equation*}
	
	If $r\geq 2$, then 
	\begin{equation*}
	P_{0}(x^r \alpha_{k}(x))=-(-1)^r \corch{e_{kk},S_{N}^{r-1}}
	-(-1)^{k+r-1}\corch{(-1)^k e_{k1},S_{N}^{k+r-2}}
		-(-1)^{N-k+r}\corch{(-1)^{N-k+1} e_{N1},S_{N}^{N-k+r-1}}
	\end{equation*}
		\begin{equation*}
-(-1)^r \paren{e_{k,k+r-1}-e_{k-r+1,k}}
	-(-1)^{r-1}\paren{e_{k,k+r-1}-0}
	-(-1)^{r+1}\paren{0- e_{k-r+1,k}}=0.
	\end{equation*}
		    If $N$ is odd 
			\begin{equation*}
		P_{0}\paren{x^r\alpha_{\frac{N+1}{2}}(x)}=-(-1)^{r}\corch{e_{\frac{N+1}{2},\frac{N+1}{2}},S_{N}^{r-1}}
			-(-1)^{\frac{N-1}{2}+r}\corch{(-1)^{\frac{N+1}{2}}\paren{e_{\frac{N+1}{2},1}+e_{N,\frac{N+1}{2}}},S_{N}^{\frac{N-3}{2}+r}}
		\end{equation*}
		\begin{equation*}
		=-(-1)^r \paren{e_{\frac{N+1}{2},\frac{N-1}{2}+r}-e_{\frac{N+3}{2}-r,\frac{N+1}{2}}}
		+(-1)^r \paren{e_{\frac{N+1}{2},\frac{N-1}{2}+r}-e_{\frac{N+3}{2}-r,\frac{N+1}{2}}}=0.
		\end{equation*}
 The second family of relations has a number of cases which we will check.

 Using Lemma~\ref{pierce}:
 \begin{itemize}
 	\item If $q<\frac{N-3}{2}$, 
 	
 	\begin{itemize}
 	\item If $2\leq k <q+2$, then $q<N-k-1$. 	Since, $k-N+q+1<0$ and $q+3<N$ we have  
 	\begin{equation*}
 		\sum_{j=0}^{q}(-1)^j S_{N}^{N+j-q-1}P_{j}(\alpha_{k})=
 	\sum_{j=0}^{k-3}(-1)^j S_{N}^{N+j-q-1}(-1)^j \paren{e_{k,k-j-1}-e_{k+j+1,k}}
 	\end{equation*}
 		\begin{equation*}
 	+(-1)^k S_{N}^{N-q+k-3}\paren{(-1)^k k e_{k1}+(-1)^{k+1}e_{2k-1}}
 	+\sum_{j=k-1}^{q}(-1)^j S_{N}^{N+j-q-1}(-1)^j \paren{e_{j+2,1}-e_{k+j+1,k}}
 	\end{equation*}
 	\begin{equation*}
 	=	\sum_{j=0}^{k-3}e_{k-N-j+q+1,k-j-1}-\sum_{j=0}^{k-3}e_{k-N+q+2,k}+
 	ke_{q+3-N,1}-e_{k-N+q+2,k}
 	\end{equation*}
 	\begin{equation*} 	
 	+\sum_{j=k-1}^{q}e_{q+3-N,1}-\sum_{j=k-1}^{q}e_{k-N+q+2,k}
 	=0.	
 	\end{equation*}

 	\item If $k=q+2$. Since, $2q-N+3<0$.
 		\begin{equation*}
 	\sum_{j=0}^{q}(-1)^j S_{N}^{N+j-q-1}P_{j}(\alpha_{k})=
 	\sum_{j=0}^{q-1}(-1)^j S_{N}^{N+j-q-1}(-1)^j \paren{e_{q+2,q+1-j}-e_{q+j+3,q+2}}
 		\end{equation*}
 	\begin{equation*}
 	+(-1)^q e_{1N}\paren{(-1)^q (q+2)e_{q+2}+(-1)^{q+1}e_{2q+3,q+2}}=
 	\sum_{j=0}^{q-1}(-1)^j S_{N}^{N+j-q-1}(-1)^j e_{2q-N-j+3,q+1-j}
 	\end{equation*}
 	 \begin{equation*}
 	-	\sum_{j=0}^{q-1}(-1)^j S_{N}^{N+j-q-1}(-1)^j e_{2q-N+4,q+2}=0.
 		\end{equation*}

 \item  If $q+2<k <N-q-1$, then $q<k+2$.  Since, $k-N+q+1<0$.
 	\begin{equation*}
 \sum_{j=0}^{q}(-1)^j S_{N}^{N+j-q-1}P_{j}(\alpha_{k})=
 \sum_{j=0}^{q}(-1)^j S_{N}^{N+j-q-1}(-1)^j \paren{e_{k,k-j-1}-e_{k+j+1,k}}
 \end{equation*}
 \begin{equation*}
 = \sum_{j=0}^{q}(-1)^j S_{N}^{N+j-q-1}(-1)^j e_{k-N-j+q+1,k-j-1}-
 \sum_{j=0}^{q}(-1)^j S_{N}^{N+j-q-1}(-1)^j e_{k-N+q+2,k}=0.
 \end{equation*}

 \item If $k=N-q-1$, then $N-k=q+1$.  Since,  $k-N+q+1<0$.
 	\begin{equation*}
 \sum_{j=0}^{q}(-1)^j S_{N}^{N+j-q-1}P_{j}(\alpha_{k})=
 \sum_{j=0}^{q-1}(-1)^j S_{N}^{N+j-q-1}P_{j}(\alpha_{k})+
 (-1)^q  S_{N}^{N-1}P_{q}(\alpha_k)
  \end{equation*}
 \begin{equation*}
 = \sum_{j=0}^{q-1}(-1)^j S_{N}^{N+j-q-1}\paren{e_{k,k-j-1}-e_{k+j+1,k}} 
   \end{equation*}
 \begin{equation*}
 +(-1)^q  e_{1N}\paren{(N-k)(-1)^{N-k+1}e_{Nk}
 +(-1)^q (e_{k,k-q-1}-e_{N,N-q-1}) } 
 \end{equation*}
 \begin{equation*}
 =\sum_{j=0}^{q-1}e_{k-N-j+q+1,k-j-1}- \sum_{j=0}^{q-1}e_{k-N+q+2,k}
 +(N-k)e_{1k}-e_{1,N-q-1}
  \end{equation*}
 \begin{equation*}
 =-(N-k-1)e_{1k}+(N-k)e_{1k}-e_{1k}=0.
 \end{equation*}

 \item  If $k>N-q-1$, then $q+1>k$.
 Since, $j\geq k-N+q+1$ implies $k-N-j+q+1\leq 0<1$.
 	\begin{equation*}
 \sum_{j=0}^{q}(-1)^j S_{N}^{N+j-q-1}P_{j}(\alpha_{k})=
 	\sum_{j=0}^{N-k-2}(-1)^j S_{N}^{N+j-q-1}(-1)^j \paren{e_{k,k-j-1}-e_{k+j+1,k}}+
 	 \end{equation*}
 	\begin{equation*}
 	(-1)^{N-k+1}  S_{N}^{2N-k-q-2}\paren{(N-k-1)(-1)^{N-k+1}e_{Nk}+(-1)^{N-k+1} e_{k,2k-N}} +
 	 	 \end{equation*}
 	\begin{equation*}
 	\sum_{j=N-k}^{q}(-1)^j S_{N}^{N+j-q-1}(-1)^j \paren{e_{k,k-j-1}-e_{N,N-j-1}}
 	=\sum_{j=0}^{N-k-2}e_{k-N-j+q+1,k-j-1}
 	-\sum_{j=0}^{N-k-2}e_{k-N+q+2,k}
 		 \end{equation*}
 	\begin{equation*}
 	+(N-k-1)e_{k-N+q+2,k}+e_{2k-2N+q+2,2k-N}
 +	\sum_{j=N-k}^{q}e_{k-N-j+q+1,k-j-1}+\sum_{j=N-k}^{q}e_{q+1-j,N-j-1}
 	 \end{equation*}
 \begin{equation*}
 =\sum_{j=0}^{q}e_{k-N-j+q+1,k-j-1}-\sum_{j=N-k}^{q}e_{q+1-j,k-j-1}=\sum_{j=q-N+k+1}^{q}e_{k-N-j+q+1,k-j-1}=0.
 \end{equation*}
 	\end{itemize}
 	
 \item  If $q=\frac{N-3}{2}$ (for $N$ odd)  then,
 \begin{itemize}
 	\item The case $2\leq k<q+2=\frac{N+1}{2}$ is similar to the case  $q<\frac{N-3}{2}$ and $2\leq k <q+2$.
 	\item If $k=q+2=\frac{N+1}{2}$ then, 
\begin{equation*}
	\sum_{j=0}^{q} (-1)^j S_{N}^{N+j-q-1}P_{j}\paren{\alpha_{\frac{N+1}{2}}}
	=\sum_{j=0}^{q-1}(-1)^j S_{N}^{N+j-q-1} (-1)^j \paren{e_{\frac{N+1}{2},\frac{N-1}{2}+j}-e_{\frac{N+3}{2}+j,\frac{N+1}{2}}}
	\end{equation*}
	\begin{equation*}
	+(-1)^{\frac{N-3}{2}}S_{N}^{N-1}(-1)^{\frac{N+1}{2}}\frac{N+1}{2}e_{\frac{N+1}{2},1}+(-1)^{\frac{N-3}{2}}S_{N}^{N-1}(-1)^{\frac{N+1}{2}}\frac{N-3}{2}e_{N,\frac{N+1}{2}}
	\end{equation*}
	\begin{equation*}
		=\sum_{j=0}^{q-1} e_{\frac{N+1}{2}-N-j+q+1,\frac{N-1}{2}-j}-\sum_{j=0}^{q-1}e_{\frac{N+3}{2}-N+q+1,\frac{N+1}{2}}
		+\frac{N-3}{2}e_{1,\frac{N+1}{2}}
	\end{equation*}
		\begin{equation*}
	=\sum_{j=0}^{q-1} e_{-j,\frac{N-1}{2}-j}-\sum_{j=0}^{q-1}e_{1,\frac{N+1}{2}}+\frac{N-3}{2}e_{1,\frac{N+1}{2}}
	=-qe_{1,\frac{N+1}{2}}+\frac{N-3}{2}e_{1,\frac{N+1}{2}}=0.
	\end{equation*}
	\item The case $q+2=\frac{N+1}{2}=N-q-1<k$ is similar to the case $q<\frac{N-3}{2}$ and $k>N-q-1$.
	\end{itemize}

	\item If $\frac{N-3}{2}< q\leq N-3$, 
	
\begin{itemize}
	\item If $0\leq k \leq N-q-2$,  then $q<N-k-1$. Since, $k-N+q+1<0$ and $q+3\leq N$. Therefore, 
	\begin{equation*}
	\sum_{j=0}^{q}(-1)^j S_{N}^{N+j-q-1}P_{j}(\alpha_{k})=
	\sum_{j=0}^{k-3}(-1)^j S_{N}^{N+j-q-1}(-1)^j \paren{e_{k,k-j-1}-e_{k+j+1,k}}
		\end{equation*}
	\begin{equation*}
	+(-1)^k S_{N}^{N-q+k-3}\paren{(-1)^k k e_{k1}+(-1)^{k+1} e_{2k-1,k} }
		+\sum_{j=k-1}^{q}(-1)^j S_{N}^{N+j-q-1}(-1)^j \paren{e_{j+2,1}-e_{k+j+1,k}}
	\end{equation*}
		\begin{equation*}
=\sum_{j=0}^{k-3} e_{k-N-j+q+1,k-j-1}- \sum_{j=0}^{k-3}	e_{k-N+q+2,k}+ke_{q+3-N,1}-e_{k-N+q+2,k}
	\end{equation*}
		\begin{equation*}
	+\sum_{j=k-1}^{q}e_{q+3-N,1}-\sum_{j=k-1}^{q}e_{k-N+q+2,k}=0.
	\end{equation*}

	\item  If $N-q-1\leq k<q+2$, then $N-k-1<q$, $k-2<q$.
	Since $q<N-k-1$ and $q+3\leq N$. Therefore,
	
		\begin{equation*}
	\sum_{j=0}^{q}(-1)^j S_{N}^{N+j-q-1}P_{j}(\alpha_{k})=
	\sum_{j=0}^{k-3}(-1)^j S_{N}^{N+j-q-1}(-1)^j \paren{e_{k,k-j-1}-e_{k+j+1,k}}
	\end{equation*}
		\begin{equation*}
	+(-1)^k S_{N}^{N-q+k-3}\paren{(-1)^k k e_{k1}+(-1)^{k+1} e_{2k-1,k} }
	+\sum_{j=k-1}^{q}(-1)^j S_{N}^{N+j-q-1}(-1)^j \paren{e_{j+2,1}-e_{k+j+1,k}}
	\end{equation*}
		\begin{equation*}
	=\sum_{j=0}^{k-3} e_{k-N-j+q+1,k-j-1}- \sum_{j=0}^{k-3}	e_{k-N+q+2,k}+ke_{q+3-N,1}-e_{k-N+q+2,k}
	+\sum_{j=k-1}^{N-k-2}e_{q+3-N,1}
		\end{equation*}
	\begin{equation*}
-\sum_{j=k-1}^{q}e_{k-N+q+2,k}	+(-1)^{N-k-1} S_{N}^{2N-k-q-2}\paren{(-1)^{N-k-1}  (N-k-1) e_{Nk}+(-1)^{N-k+1} e_{N-k+1,1} }
		\end{equation*}
	\begin{equation*}
	+\sum_{j=N-k}^{q}(-1)^j S_{N}^{N+j-q-1}(-1)^j \paren{e_{N,N-j-1}-e_{j+2,1}}
	\end{equation*}
		\begin{equation*}
	=\sum_{j=0}^{k-3} e_{k-N-j+q+1,k-j-1}- \sum_{j=0}^{k-3}	e_{k-N+q+2,k}+ke_{q+3-N,1}-e_{k-N+q+2,k}
	+\sum_{j=k-1}^{N-k-2}e_{q+3-N,1}
	\end{equation*}
	\begin{equation*}
-\sum_{j=k-1}^{q}e_{k-N+q+2,k}	+  (N-k-1) e_{k-N+q+2,k}+ e_{q+3-N,1} +\sum_{j=N-k}^{q}e_{q+1-j,N-j-1}-\sum_{j=N-k}^{q}e_{q+3-N,1}
	\end{equation*}
		\begin{equation*}
	=\sum_{j=0}^{k-3} e_{k-N-j+q+1,k-j-1}+\sum_{j=N-k}^{q}e_{q+1-j,N-j-1}=0.
	\end{equation*}

	\item If $k=N-q-1$, then $N-k-1=q$.
	Since $q+3\leq N$.
		\begin{equation*}
	\sum_{j=0}^{q}(-1)^j S_{N}^{N+j-q-1}P_{j}(\alpha_{k})=
	\sum_{j=0}^{k-3}(-1)^j S_{N}^{N+j-q-1}(-1)^j \paren{e_{k,k-j-1}-e_{k+j+1,k}}
	\end{equation*}
	\begin{equation*}
	+(-1)^k S_{N}^{N-q+k-3}\paren{(-1)^k k e_{k1}+(-1)^{k+1} e_{2k-1,k} }
	+\sum_{j=k-1}^{N-k-2}(-1)^j S_{N}^{N+j-q-1}(-1)^j \paren{e_{j+2,1}-e_{k+j+1,k}}
		\end{equation*}
	\begin{equation*}
	+(-1)^q S_{N}^{N-1}\paren{(N-k-1)(-1)^{N-k-1}e_{Nk}+(-1)^{N-k+1}e_{N-k+1,1}}
	\end{equation*}
		\begin{equation*}
	=\sum_{j=0}^{k-3} e_{k-N-j+q+1,k-j-1}- \sum_{j=0}^{k-3}	e_{k-N+q+2,k}+ke_{q+3-N,1}-e_{k-N+q+2,k}
	 	\end{equation*}
 	\begin{equation*}
	+\sum_{j=k-1}^{N-k-2}e_{q+3-N,1}-\sum_{j=k-1}^{N-k-2}e_{k-N+q+2,k}	+ (N-k-1) e_{1k}
	\end{equation*}
	\begin{equation*}
	=\sum_{j=0}^{k-3} e_{-j,k-j-1}- \sum_{j=0}^{k-3}	e_{1k}-e_{1k}
	-\sum_{j=k-1}^{N-k-2}e_{1k}	+ (N-k-1) e_{1k}=0.
	\end{equation*}

	\item If $k=q+2$, then $k-2=q$. Since $q<N-k-1$.
	\begin{equation*}
	\sum_{j=0}^{q}(-1)^j S_{N}^{N+j-q-1}P_{j}(\alpha_{k})=
	\sum_{j=0}^{k-3}(-1)^j S_{N}^{N+j-q-1}(-1)^j \paren{e_{k,k-j-1}-e_{k+j+1,k}}
	\end{equation*}
	\begin{equation*}
	+(-1)^k S_{N}^{N-q+k-3}\paren{(-1)^k k e_{k1}+(-1)^{k+1} e_{2k-1,k} }
	=\sum_{j=0}^{k-3} e_{k-N-j+q+1,k-j-1}- \sum_{j=0}^{k-3}	e_{k-N+q+2,k}
		\end{equation*}
	\begin{equation*}
	+ke_{q+3-N,1}-e_{k-N+q+2,k}=0.
	\end{equation*}

	\item If $k>q+2$, then $k>q>N-q-1$ and $N-k-1<q$
		\begin{equation*}
	\sum_{j=0}^{q}(-1)^j S_{N}^{N+j-q-1}P_{j}(\alpha_{k})=
	\sum_{j=0}^{N-k-2}(-1)^j S_{N}^{N+j-q-1}(-1)^{j} \paren{e_{k,k-j-1}-e_{k+j+1,k}}
	\end{equation*}
	\begin{equation*}
+(-1)^{N-k-1} S_{N}^{2N-k-q-2}\paren{(-1)^{N-k+1} (N-k) e_{Nk}+(-1)^{N-k+1} \paren{e_{k,2k-N}-e_{Nk}} }
	\end{equation*}
\begin{equation*}
+\sum_{j=N-k}^{q}(-1)^j S_{N}^{N+j-q-1}(-1)^{j} \paren{e_{k,k-j-1}-e_{N,N-j-1}}
\end{equation*}
\begin{equation*}
=\sum_{j=0}^{N-k-2}e_{k-N-j+q+1,k-j-1}-\sum_{j=0}^{N-k-2}e_{k-N+q+2,k}
+ (N-k) e_{k-N+q+2,k}+e_{2k-2N+q+2,2k-N}
\end{equation*}
\begin{equation*}
-e_{k-N+q+2,k}+\sum_{j=N-k}^{q}e_{k-N-j+q+1,k-j-1}-\sum_{j=N-k}^{q}e_{q+1-j,N-j-1}
=\sum_{j=0}^{q}e_{k-N-j+q+1,k-j-1}
\end{equation*}
\begin{equation*}
-\sum_{j=N-k}^{q}e_{q+1-j,N-j-1}=\sum_{j=k-N+q+1}^{q}e_{k-N-j+q+1,k-j-1}=0.
\end{equation*}
 \end{itemize}
 \item If $q=N-2$ and $k<\frac{N+1}{2}$ then $k-1<N-k\leq N-1$. Therefore, 
 \begin{equation*}
     \sum_{j=0}^{N-2}(-1)^j S_{N}^{N+j-q-1}P_{j}(\alpha_{k})
     =\sum_{j=0}^{N-2}(-1)^j S_{N}^{j+1}P_{j}(\alpha_{k})
     =\sum_{j=0}^{k-3}(-1)^j S_{N}^{j+1}(-1)^j (e_{k,k-j-1}-e_{k+j+1,k})
     \end{equation*}
\begin{equation*}
     +(-1)^k S_{N}^{k-1}\llave{(-1)^k k e_{k1}+(-1)^{k+1}e_{2k-1,k}}
     +\sum_{j=k-1}^{N-k-2}(-1)^j S_{N}^{j+1}(-1)^j (e_{j+2,1}-e_{k+j+1,k})
 \end{equation*}
 \begin{equation*}
    +(-1)^{N-k+1}S_{N}^{N-k}\llave{(N-k-1)(-1)^{N-k-1}e_{Nk}+(-1)^{N-k+1}e_{N-k+1,1}}
     \end{equation*}
 \begin{equation*}
    +\sum_{j=N-k}^{N-3}(-1)^j S_{N}^{j+1}(-1)^{j+1} (e_{N,N-j-1}-e_{j+2,1})
    +(-1)^{N}S_{N}^{N-1} (N-1) (-1)^{N-1} e_{N1}
 \end{equation*}
\begin{equation*}
    =\sum_{j=0}^{k-3} e_{k-j-1,k-j-1} -\sum_{j=0}^{k-3} e_{k,k}+ ke_{11} -e_{kk} + \sum_{j=k-1}^{N-k-2} e_{11}
    -\sum_{j=k-1}^{N-k-2} e_{kk} +(N-k-1) e_{kk}
 \end{equation*}
\begin{equation*}    
    +e_{11}-\sum_{j=N-k}^{N-3}e_{N-j-1,N-j-1}+\sum_{j=N-k}^{N-3}e_{11}
    -(N-1)e_{11}
\end{equation*}
\begin{equation*}
= \sum_{j=0}^{k-3} e_{k-j-1,k-j-1}- (k-2)e_{kk}+ke_{11}-e_{kk}+ (N-k)e_{11}-(N-2k)e_{kk}
+e_{11}
\end{equation*}
\begin{equation*}
-\sum_{r=0}^{k-3} e_{k-r-1,k-r-1}+(k-2)e_{11}-(N-1)e_{11}=0.
\end{equation*}

\item If $q=N-1$ and $k< \frac{N+1}{2}$ then 
\begin{equation*}
     \sum_{j=0}^{N-1}(-1)^j S_{N}^{N+j-q-1}P_{j}(\alpha_{k})
     =\sum_{j=0}^{N-1}(-1)^j S_{N}^{j}P_{j}(\alpha_{k})
     =\sum_{j=0}^{k-3} (-1)^{j} S_{N}^{j} (-1)^{j} (e_{k,k-j-1}-e_{k+j+1,k})
     \end{equation*}
\begin{equation*}
     +(-1)^k S_{N}^{k-2}\llave{(-1)^k k e_{k1}+ (-1)^{k+1}e_{2k-1,k}}
     +\sum_{j=k-1}^{N-k-2} (-1)^{j} S_{N}^{j} (-1)^{j} (e_{j+2,1}-e_{k+j+1,k})
     \end{equation*}
\begin{equation*}
     +(-1)^{N-k-1}S_{N}^{N-k-1}\llave{(N-k-1)(-1)^{N-k-1}e_{Nk}+(-1)^{N-k+1}e_{N-k+1,1}}
\end{equation*}
\begin{equation*}
+\sum_{j=N-k}^{N-3} (-1)^{j} S_{N}^{j} (-1)^{j+1} (e_{N,N-j-1}-e_{j+2,1})
+(-1)^{N-2}S_{N}^{N-2}(N-1)(-1)^{N-1}e_{N1}
\end{equation*}
\begin{equation*}
    =\sum_{j=0}^{k-3} e_{k-j,k-j-1}-\sum_{j=0}^{k-3} e_{k+1,k}
    +ke_{21}-e_{k+1,k}+\sum_{j=k-1}^{N-k-2}e_{21}
    +\sum_{j=k-1}^{N-k-2}e_{k+1,k}
    +(N-k-1)e_{k+1,k}
    \end{equation*}
\begin{equation*}
    +e_{21}-\sum_{j=N-k}^{N-3}e_{N-j,N-j-1}+ \sum_{j=N-k}^{N-3}e_{21}-(N-1)e_{21}
\end{equation*}
\begin{equation*}
=\sum_{j=0}^{k-3} e_{k-j,k-j-1}-(k-2)e_{k+1,k}+ke_{21}-e_{k+1,k}
+(N-2k)e_{21}-(N-2k)e_{k+1,k}+(N-k-1)e_{k+1,k}
\end{equation*}
\begin{equation*}
+e_{21}
-\sum_{r=0}^{k-3} e_{k-r,k-r-1}+(k-2)e_{21}-(N-1)e_{21}=0.
\end{equation*}
\item If $q=N-2$ and $K>\frac{N+1}{2}$ then $N-k<k-1\leq N-1$ and 
\begin{equation*}
        \sum_{j=0}^{N-2}(-1)^j S_{N}^{N+j-q-1}P_{j}(\alpha_{k})
        =     \sum_{j=0}^{N-2}(-1)^j S_{N}^{j+1}P_{j}(\alpha_{k})
        =  \sum_{j=0}^{N-k-2}(-1)^j S_{N}^{j+1} (-1)^j (e_{k,k-j-1}-e_{k+j+1,k})
 \end{equation*}
\begin{equation*}       
        +(-1)^{N-k-1}S_{N}^{N-k}\llave{(N-k-1)(-1)^{N-k-1}e_{Nk}+(-1)^{N-k+1}e_{k,2k-N}}
\end{equation*}
\begin{equation*}
+\sum_{j=N-k}^{k-3}(-1)^j S_{N}^{j+1} (-1)^j (e_{k,k-j-1}-e_{N,N-j-1})
+(-1)^{k}S_{N}^{k-1}\llave{k(-1)^ke_{k1}+(-1)^{k+1}e_{N,N-k+1}}
\end{equation*}
\begin{equation*}
+\sum_{j=k-1}^{N-3}(-1)^j S_{N}^{j+1} (-1)^{j+1} (e_{N,N-j-1}-e_{j+2,1})
+(-1)^{N}S_{N}^{N-1}(N-1)(-1)^{N-1}e_{N1}.
\end{equation*}
\begin{equation*}
= \sum_{j=0}^{N-k-2}e_{k-j-1,k-j-1}- \sum_{j=0}^{N-k-2}e_{kk}+(N-k-1)e_{kk}
+e_{2k-n,2k-N}+ \sum_{j=N-k}^{k-3}e_{k-j-1,k-j-1}
\end{equation*}
\begin{equation*}
- \sum_{j=N-k}^{k-3}e_{N-j-1,N-j-1}
+ke_{11}-e_{N-k-1,N-k-1}- \sum_{j=k-1}^{N-3}e_{N-j-1,N-j-1}+ \sum_{j=k-1}^{N-3}e_{11}
-(N-1)e_{11}
\end{equation*}
\begin{equation*}
    =\sum_{j=0}^{k-3}e_{k-j-1,k-j-1}-(N-k-1)e_{kk}+(N-k-1)e_{kk}
    -\sum_{j=N-k}^{N-3}e_{N-j-1,N-j-1}+ke_{11}
\end{equation*}
\begin{equation*}
+(N-k-1)e_{11}-(N-1)e_{11}= \sum_{j=0}^{k-3}e_{k-j-1,k-j-1}-\sum_{r=0}^{k-3}e_{k-r-1,k-r-1}=0.
\end{equation*}

\item If $q=N-1$ and $k>\frac{N+1}{2}$ then $N-k<k-1\leq N-1$ and 

\begin{equation*}
     \sum_{j=0}^{N-1}(-1)^j S_{N}^{N+j-q-1}P_{j}(\alpha_{k})
     = \sum_{j=0}^{N-1}(-1)^j S_{N}^{j}P_{j}(\alpha_{k})
     =\sum_{j=0}^{N-k-2}(-1)^j S_{N}^{j} (-1)^{j} (e_{k,k-j-1}-e_{k+j+1,k})
      \end{equation*}
\begin{equation*}  
     + (-1)^{N-k-1}S_{N}^{N-k-1}\llave{(N-k-1)(-1)^{N-k-1}e_{Nk}+(-1)^{N-k-1}e_{k,2k-N}}
        \end{equation*}
\begin{equation*}     
      +\sum_{j=N-k}^{k-3}(-1)^j S_{N}^{j} (-1)^{j} (e_{k,k-j-1}-e_{N,N-j-1})
       \end{equation*}
\begin{equation*}  
     +(-1)^{k}S_{N}^{k-2}\llave{k(-1)^{k}e_{k1}+(-1)^{k+1}e_{N,N-k+1}}
     +\sum_{j=k-1}^{N-3}(-1)^j S_{N}^{j} (-1)^{j+1} (e_{N,N-j-1}-e_{j+2,1})
      \end{equation*}
\begin{equation*}  
    +(-1)^{N}S_{N}^{N-2}(N-1)(-1)^{N-1}e_{N1}
=\sum_{j=0}^{N-k-2} e_{k-j, k-j-1}- \sum_{j=0}^{N-k-2} e_{k+1,k} + (N-k-1) e_{k+1,k}
 \end{equation*}
\begin{equation*}  
+e_{2k-N+1, 2k-N} +\sum_{j=N-k}^{k-3}e_{k-j,k-j-1}-\sum_{j=N-k}^{k-3}e_{N-j,N-j-1}
+ke_{21}-e_{N-k+2,N-k+1}-\sum_{j=k-1}^{N-3}e_{N-j,N-j-1}
\end{equation*}
\begin{equation*}
+\sum_{j=k-1}^{N-3}e_{21}-(N-1)e_{21}= \sum_{j=0}^{k-3}e_{k-j,k-j-1}-(N-k-1)e_{k+1,k}+(N-k-1)e_{k-1,k}
-\sum_{j=N-k}^{N-3}e_{N-j,N-j-1}
\end{equation*}
\begin{equation*}
+ke_{21}+(N-k-1)e_{21}-(N-1)e_{21}
=\sum_{j=0}^{k-3}e_{k-j,k-j-1}-\sum_{r=0}^{k-3}e_{k-r,k-r-1}=0.
\end{equation*}

\item If $q=N-2$ and $k=\frac{N+1}{2}$ for $N$ odd we have $N-k=k-1<N-1$, then

\begin{equation*}
  \sum_{j=0}^{N-2}(-1)^j S_{N}^{N+j-q-1}P_{j}(\alpha_{k})
  = \sum_{j=0}^{N-2}(-1)^j S_{N}^{j+1}P_{j}(\alpha_{k})
  = \sum_{j=0}^{\frac{N-5}{2}}(-1)^j S_{N}^{j+1} (-1)^{j} (e_{\frac{N+1}{2},\frac{N-1}{2}-j}-e_{\frac{N+3}{2}+j,\frac{N+1}{2}})
\end{equation*}
\begin{equation*}  
  +(-1)^{\frac{N-3}{2}}S_{N}^{\frac{N-1}{2}}\llave{\paren{\frac{N-3}{2}}(-1)^{\frac{N+1}{2}}e_{N,\frac{N+1}{2}}+(-1)^{\frac{N+1}{2}}\paren{\frac{N+1}{2}}
  e_{\frac{N+1}{2},1}}
\end{equation*}
\begin{equation*}
    +\sum_{j=\frac{N-1}{2}}^{N-3}(-1)^j S_{N}^{j+1} (-1)^{j+1}(e_{N,N-j-1}-e_{j+2,1})
    +(-1)^{N-2}S_{N}^{N-1}(N-1)(-1)^{N+1}e_{N1}
\end{equation*}
\begin{equation*}
    =\sum_{j=0}^{\frac{N-5}{2}}e_{\frac{N-1}{2}-j,\frac{N-1}{2}-j}- \sum_{j=0}^{\frac{N-5}{2}}e_{\frac{N+1}{2},\frac{N+1}{2}}
    +\paren{\frac{N-3}{2}}e_{\frac{N+1}{2},\frac{N+1}{2}}+ \paren{\frac{N+1}{2}}e_{\frac{N+1}{2},1}
    -\sum_{j=\frac{N-1}{2}}^{N-3} e_{N-j-1,N-j-1}
 \end{equation*}
\begin{equation*}   
+ \sum_{j=\frac{N-1}{2}}^{N-3}e_{11}    -(N-1)e_{11}
    = \sum_{j=0}^{\frac{N-5}{2}}e_{\frac{N-1}{2}-j,\frac{N-1}{2}-j}-\paren{\frac{N-3}{2}}e_{\frac{N+1}{2},\frac{N+1}{2}}
    +\paren{\frac{N-3}{2}}e_{\frac{N+1}{2},\frac{N+1}{2}}
 \end{equation*}
\begin{equation*}    
+\paren{\frac{N+1}{2}}e_{11}    -\sum_{r=0}^{\frac{N-5}{2}}e_{\frac{N-1}{2}-r,\frac{N-1}{2}-r}
    +\paren{\frac{N-3}{2}}e_{11}-(N-1)e_{11}=0.
\end{equation*}

\item If $q=N-1$ and $k=\frac{N+1}{2}$ then $N-k=k-1<N-1$ and 
\begin{equation*}
     \sum_{j=0}^{N-1}(-1)^j S_{N}^{N+j-q-1}P_{j}(\alpha_{k})
     = \sum_{j=0}^{N-1}(-1)^j S_{N}^{j}P_{j}(\alpha_{k})
     =\sum_{j=0}^{\frac{N-5}{2}}(-1)^{j}S_{N}^j (-1)^{j} (e_{\frac{N+1}{2},\frac{N-1}{2}-j}-e_{\frac{N+3}{2}+j,\frac{N+1}{2}})
  \end{equation*}
\begin{equation*}       
     +(-1)^\frac{N-3}{2} S_{N}^{\frac{N-3}{2}}\llave{\paren{\frac{N-3}{2}}(-1)^{\frac{N+1}{2}}e_{N,\frac{N+1}{2}}
     +(-1)^{\frac{N+1}{2}}\paren{\frac{N+1}{2}}e_{\frac{N+1}{2},1}}
   \end{equation*}
\begin{equation*}        
     +\sum_{j=\frac{N-1}{2}}^{N-3} (-1)^{j}S_{N}^{j}(-1)^{j+1}(e_{N,N-j-1}-e_{j+2,1})
     +(-1)^{N-2}S_{N}^{N-2}(N-1)(-1)^{N+1}e_{N1}
       \end{equation*}
\begin{equation*}  
     =\sum_{j=0}^{\frac{N-5}{2}} e_{\frac{N+1}{2}-j,\frac{N-1}{2}-j}-
     \sum_{j=0}^{\frac{N-5}{2}}e_{\frac{N+3}{2},\frac{N+1}{2}}+\paren{\frac{N-3}{2}}e_{\frac{N+3}{2},\frac{N+1}{2}}
    + \paren{\frac{N+1}{2}}e_{21}-\sum_{j=\frac{N-1}{2}}^{N-3} e_{N-j,N-j-1} 
       \end{equation*}
\begin{equation*}   
    +\sum_{j=\frac{N-1}{2}}^{N-3}e_{21}-(N-1)e_{21} 
    =\sum_{j=0}^{\frac{N-5}{2}} e_{\frac{N+1}{2}-j,\frac{N-1}{2}-j}
    -\paren{\frac{N-3}{2}}e_{\frac{N+3}{2},\frac{N+1}{2}}+\paren{\frac{N-3}{2}}e_{\frac{N+3}{2},\frac{N+1}{2}}
        \end{equation*}
\begin{equation*}     
 + \paren{\frac{N+1}{2}}e_{21}
    -\sum_{r=0}^{\frac{N-5}{2}} e_{\frac{N+1}{2}-r,\frac{N-1}{2}-r}
    + \paren{\frac{N-3}{2}}e_{21}-(N-1)e_{21}=0. 
\end{equation*}

 \end{itemize}

\end{demos}


Recall that a Pierce decomposition of a noncommutative ring $R$ with unit $1$ is a finite set of elements $r_1, ... ,r_n \in R$ such that  $1=\sum_{j=1}^{n}r_j$ and $r_i r_j =\delta_{ij}$ for all $1\leq i,j \leq n$. Now we state the Pierce decomposition of $\mathbb{A}$.
\begin{coro}
	If $\alpha_{1}=I-\sum_{k=2}^{N-1}\alpha_{k}$, then $\llave{\alpha_k}_{1\leq k \leq N-1}$ is a Pierce decomposition of $\mathbb{A}$.
\end{coro}
\begin{demos}
	In fact, if $2\leq k<l\leq N-1$, then 
	\begin{equation*}
	\alpha_{k}\alpha_l =\paren{e_{kk}+\sum_{j=1}^{N-1}a_{kj}x^j}\paren{e_{ll}+\sum_{r=1}^{N-1}a_{lr}x^r}
	=\sum_{j,r=1}^{N-1}a_{kj}a_{lr}x^{j+r}.
	\end{equation*}
	However, 
	{\footnotesize
	\begin{equation*}
	  	a_{kj}a_{lr}=\paren{(-1)^{j+1}(\delta_{k,j+1}e_{k1}+\delta_{k,N-j}e_{Nk})+(-1)^j \delta_{N,j+1}e_{N1}}
	\paren{(-1)^{r+1}(\delta_{l,r+1}e_{l1}+\delta_{l,N-r}e_{Nl})+(-1)^j \delta_{N,r+1}e_{N1}}
	=0
	\end{equation*}}
 implies $\alpha_{k}\alpha_l =0$.
	
	On the other hand, 
	\begin{equation*}
	\alpha_{k}^2=\paren{e_{kk}+\sum_{j=1}^{N-1}a_{kj}x^j}\paren{e_{kk}+\sum_{j=1}^{N-1}a_{kj}x^j}
	=e_{kk}+\sum_{j=1}^{N-1}(e_{kk}a_{kj}+a_{kj}e_{kk})x^j + \sum_{j,r=1}^{N-1}a_{kj}a_{kr}x^{j+r}.
	\end{equation*}
	However, $$e_{kk}a_{kj}+a_{kj}e_{kk}=	(-1)^{j+1}\paren{\delta_{k,j+1}e_{k1}+\delta_{k,N-j}e_{Nk}}$$ and 
	$$a_{kj}a_{kr}=(-1)^{j+r}\delta_{k,N-j}\delta_{k,r+1}e_{N1}
	=(-1)^{N-1}\delta_{j,N-k}\delta_{r,N-k}e_{N1}$$
	imply 
		\begin{equation*}
	\alpha_{k}^2=e_{kk}+\sum_{j=1}^{N-1} (-1)^{j+1}\paren{\delta_{k,j+1}e_{k1}+\delta_{k,N-j}e_{Nk}}x^j+(-1)^{N-1}e_{N1}x^{N-1}
	=e_{kk}+\sum_{j=1}^{N-1}a_{kj}x^j=\alpha_k.
	\end{equation*}
	Thus, $\alpha_k \alpha_l =\delta_{kl}\alpha_k$ for $2\leq k,l \leq N-1$. By the definition of $\alpha_1$ and the previous properties we can extend this to 
	$\alpha_k \alpha_l =\delta_{kl}\alpha_k$ for $1\leq k,l \leq N-1$. The assertion follows by the Theorem \ref{pierce A}.
\end{demos}

Corollaries \ref{first1}, \ref{second}, and Proposition~\ref{calogero} give an answer to the Conjectures 1, 2 and 3  of \cite{Grfrm-4} about three bispectral full rank 1 algebras as we shall describe in the following sections. Moreover, these algebras are Noetherian and finitely generated because they are contained in the $N\times N$ matrix polynomial ring $M_{N}(\mathbb{K}[x])$.
\begin{coro}\label{first1}
	Let $\Gamma$ be the sub-algebra of $M_{2}(\complex)[x]$ of the form 
	\begin{equation*}
	\paren{
		\begin{matrix} 
		r_{0}^{11} & r_{0}^{12} \\
		0 & r_{0}^{11} 
		\end{matrix}}+  \paren{
		\begin{matrix} 
		r_{1}^{11} & r_{1}^{12} \\
		0 & r_{1}^{11} 
		\end{matrix}}x+ \paren{
		\begin{matrix} 
		r_{2}^{11} & r_{2}^{12} \\
		r_{1}^{11} & r_{2}^{22} 
		\end{matrix}}x^{2}+
	\paren{
		\begin{matrix} 
		r_{3}^{11} & r_{3}^{12} \\
		r_{2}^{22}+r_{2}^{11}-r_{1}^{12} & r_{3}^{22} 
		\end{matrix}}x^3+x^{4}p(x),
	\end{equation*} 
	where $p\in M_{2}(\complex)[x]$ and all the variables $r_{0}^{11},r_{0}^{12}, r_{1}^{11},r_{1}^{12},r_{2}^{11},r_{2}^{22},r_{3}^{11},r_{3}^{12}, r_{3}^{22} \in \complex$. Then $\Gamma=\mathbb{A}$. Moreover, for each $\theta$ we have an explicit expression for the operator $\mathcal{B}$.
\end{coro}
\begin{demos}
The proof is given by the Theorem \ref{Gen} with $N=2$. \qed
\end{demos}


\begin{coro}\label{second}
	Let $\Gamma$ the sub-algebra of $M_{3}(\complex)[x]$ of the form 
	\begin{equation*}
	\paren{
		\begin{matrix} 
		r_{0}^{11} & r_{0}^{12} & r_{0}^{13} \\
		0 & r_{0}^{22} & r_{0}^{23} \\
		0 & 0 & r_{0}^{11}
		\end{matrix}}+ \paren{
		\begin{matrix} 
		r_{1}^{11} & r_{1}^{12} & r_{1}^{13} \\
		r_{0}^{22}-r_{0}^{11} & r_{1}^{22} & r_{1}^{23} \\
		0 & r_{0}^{22}-r_{0}^{11} & r_{1}^{11}+r_{0}^{23}-r_{0}^{12}
		\end{matrix}}x
	\end{equation*}
	\begin{equation*}
	+ \paren{
		\begin{matrix} 
		r_{2}^{11} & r_{2}^{12} & r_{2}^{13} \\
		r_{1}^{22}-r_{1}^{11}-r_{0}^{23}+r_{0}^{12} & r_{2}^{22} & r_{2}^{23} \\
		r_{0}^{22}-r_{0}^{11} & r_{1}^{22}-r_{1}^{11} & r_{2}^{11}+r_{1}^{23}-r_{1}^{12}
		\end{matrix}}x^{2}
	+  \paren{
		\begin{matrix} 
		r_{3}^{11} & r_{3}^{12} & r_{3}^{13} \\
		r_{3}^{21}& r_{3}^{22} & r_{3}^{23} \\
		r_{1}^{22}-2r_{1}^{11}-r_{0}^{23}+r_{0}^{12} & r_{3}^{32} & r_{3}^{33}
		\end{matrix}}x^{3}   
	\end{equation*}
	\begin{equation*}
	+ \paren{
		\begin{matrix} 
		r_{4}^{11} & r_{4}^{12} & r_{4}^{13} \\
		r_{4}^{21} & r_{4}^{22} & r_{4}^{23} \\
		r_{3}^{32}+r_{3}^{21}-r_{2}^{22}-r_{2}^{11}+r_{1}^{12} & r_{4}^{22} & r_{4}^{33}
		\end{matrix}}x^{4}
	\end{equation*}	
	\begin{equation*}
+  \paren{
	\begin{matrix} 
	r_{5}^{11} & r_{5}^{12} & r_{5}^{13} \\
	r_{5}^{21}& r_{5}^{22} & r_{5}^{23} \\
	r_{4}^{32}+r_{4}^{21}-r_{3}^{33}-r_{3}^{22}-r_{3}^{11}+r_{2}^{23}+r_{2}^{12}-r_{1}^{13} & r_{5}^{32} & r_{5}^{33}
	\end{matrix}}x^{5}	+x^{6}p(x) \mbox{ , }
	\end{equation*} \label{set}
	where $p\in M_{3}(\complex)[x]$ and all the variables $r_{0}^{11}, r_{0}^{12},...,r_{5}^{33}\in \complex$ are arbitrary.
	
	Then, $\Gamma =\mathbb{A}$ and for each $\theta$ we have  an explicit expression for the operator $\mathcal{B}$.
\end{coro}
\begin{demos}
The proof is given by the Theorem \ref{Gen} with $N=3$. \qed
\end{demos}

\subsection{An Example linked to the Spin Calogero Systems}\label{spin}

This example is linked to the spin Calogero  systems whose relation with bispectrality can be found in \cite{BGK09}. We consider the case when both "eigenvalues" $F$ and $\theta$ are matrix valued. 
Let 
\begin{equation*}
\psi(x,z)=\frac{e^{xz}}{(x-2)xz}\paren{\begin{matrix} 
\frac{x^3z^2-2x^2z^2-2x^2z+3xz+2x-2}{xz} & \frac{1}{x} \\
\frac{xz-2}{z} & x^2z-2xz-x+1
\end{matrix}}
\end{equation*}
and \begin{equation*}
\mathcal{L}=\paren{\begin{matrix} 
     0 & 0 \\
	0 & 1
	\end{matrix}}.\partial_{x}^2
+\paren{\begin{matrix} 
0 & \frac{1}{(x-2)x^2} \\
	-\frac{1}{x-2} & 0
	\end{matrix}}.\partial_{x}
+ \paren{\begin{matrix} 
	-\frac{1}{x^2(x-2)^2} & \frac{x-1}{x^3(x-2)^2} \\
	\frac{2x-1}{x(x-2)^2} & -\frac{2x^2-4x+3}{x^2(x-2)^2}
	\end{matrix}},
\end{equation*}
then $\mathcal{L}\psi=\psi F$ with 
\begin{equation*}
F(z)=\paren{\begin{matrix} 
	0 & 0 \\
	0 & z^2
	\end{matrix}}.
\end{equation*}
On the other hand, it  is easy to check that $\psi \mathcal{B}=\theta \psi$ for 
\begin{equation*}
\mathcal{B}=\partial_{z}^3.\paren{\begin{matrix} 
                   0 & 0 \\
	               1 & 0
	\end{matrix}}+ \partial_{z}^2.
\paren{\begin{matrix} 
	0 & 0 \\
	-\frac{2z+1}{z} & 0
	\end{matrix}}+
	\partial_{z}.\paren{\begin{matrix} 
		1 & 0 \\
		\frac{2(z-1)}{z^2} & 1
		\end{matrix}}
		+\paren{\begin{matrix} 
			-z^{-1} & 0 \\
			6z^{-3} & z^{-1}
			\end{matrix}}
\end{equation*}
and 
\begin{equation*}
\theta(x)=
\paren{\begin{matrix} 
	x & 0 \\
	x^2(x-2) & x
	\end{matrix}}.
\end{equation*} 
The following proposition characterizes the algebra $\mathbb{A}$ of all polynomial $F$ such that there exist $\mathcal{L}=\mathcal{L}(x,\partial_{x})$ with $\mathcal{L}\psi=\psi F$.
\begin{teor}\label{calogero}
	Let $\Gamma$ be the sub-algebra of $M_{2}(\complex)[z]$ of the form 
	\begin{equation*}
	\paren{\begin{matrix} 
		a & 0 \\
		b-a & b
		\end{matrix}}+
	\paren{\begin{matrix} 
		c & c \\
		a-b-c & -c
		\end{matrix}}z+
	\paren{\begin{matrix} 
		a-b-c & c+a-b\\
		d & e
		\end{matrix}}\frac{z^2}{2}+z^3p(z),
	\end{equation*}
	where $p\in M_{2}(\complex)[z]$ and all the variables $a,b,c,d,e$ are arbitrary. Then $\Gamma=\mathbb{A}$.
	
	\end{teor}
	
\begin{demos}
	We shall break the proof in different steps. 

Step 1:	The set $\Gamma$ is an algebra.
Clearly if $F_{1}, F_{2}\in \Gamma$, then $F_{1}+F_{2}\in \Gamma$ and $\alpha F_{1}\in \Gamma$ if $\alpha \in \complex$,
\begin{equation*}
F_{1}(z)=
\paren{\begin{matrix} 
	a_{1} & 0 \\
	b_{1}-a_{1} & b_{1}
	\end{matrix}}+
\paren{\begin{matrix} 
	c_{1} & c_{1} \\
	a_{1}-b_{1}-c_{1} & -c_{1}
	\end{matrix}}z+
\paren{\begin{matrix} 
	a_{1}-b_{1}-c_{1} & c_{1}+a_{1}-b_{1}\\
	d_{1} & e_{1}
	\end{matrix}}\frac{z^2}{2}+z^3p_{1}(z),
\end{equation*}
\begin{equation*}
F_{2}(z)=
\paren{\begin{matrix} 
	a_{2} & 0 \\
	b_{2}-a_{2} & b_{2}
	\end{matrix}}+
\paren{\begin{matrix} 
	c_{2} & c_{2} \\
	a_{2}-b_{2}-c_{2} & -c_{2}
	\end{matrix}}z+
\paren{\begin{matrix} 
	a_{2}-b_{2}-c_{2} & c_{2}+a_{2}-b_{2}\\
	d_{2} & e_{2}
	\end{matrix}}\frac{z^2}{2}+z^3p_{2}(z).
\end{equation*}

Thus,
\begin{equation*}
F_{1}(z)F_{2}(z)=\paren{\begin{matrix} 
	a_{1}a_{2} & 0 \\
	b_{1}b_{2}-a_{1}a_{2} & b_{1}b_{2}
	\end{matrix}}+
\paren{\begin{matrix} 
	a_{1}c_{2}+b_{2}c_{1} & a_{1}c_{2}+b_{2}c_{1} \\
	a_{1}a_{2}-b_{1}b_{2}-a_{1}c_{2}-b_{2}c_{1} & -(a_{1}c_{2}+b_{2}c_{1})
	\end{matrix}}z+
\end{equation*}
\begin{equation*}
\left(\begin{matrix} 
a_{1}a_{2}-b_{1}b_{2}-a_{1}c_{2}-b_{2}c_{1} \\
b_{2}e_{1}-a_{2}e_{1}+b_{1}d_{2}+a_{2}d_{1}-3b_{1}c_{2}+3a_{1}c_{2}+2b_{2}c_{1}-2a_{2}c_{1}-b_{1}b_{2}+a_{1}b_{2}+a_{2}b_{1}-a_{1}a_{2} 
\end{matrix}\right.
\end{equation*}
\begin{equation*}
\left. \begin{matrix} 
a_{1}c_{2}+b_{2}c_{1}+a_{1}a_{2}-b_{1}b_{2}\\
b_{1}e_{2}+b_{2}e_{1}-b_{1}c_{2}+a_{1}c_{2}-b_{1}b_{2}+a_{1}b_{2}+a_{2}b_{1}-a_{1}a_{2}
\end{matrix}\right)\frac{z^2}{2} +z^3p(z).
\end{equation*}
For some polynomial $p\in M_{2}(\complex)[z]$. In particular $F_{1}F_{2}\in \Gamma$. Since $M_{2}(\complex)[z]$ is an algebra and $\Gamma$ is closed for operations induced by $M_{2}(\complex)[z]$ we have that $\Gamma$ is an algebra.

Step 2: There exists a finite dimensional vector space $E$ such that $\Gamma =E \oplus z^3 M_{2}(\complex)[z]$.

Consider $F\in \Gamma$, then 
\begin{equation*}
F(z)=\paren{\begin{matrix} 
	a & 0 \\
	b-a & b
	\end{matrix}}+
\paren{\begin{matrix} 
	c & c \\
	a-b-c & -c
	\end{matrix}}z+
\paren{\begin{matrix} 
	a-b-c & c+a-b\\
	d & e
	\end{matrix}}\frac{z^2}{2}+z^3p(z)
\end{equation*}
\begin{equation*}
=a\alpha_{1}+b\alpha_{2}+c\alpha_{3}+d\alpha_{4}+e\alpha_{5}+z^3p(z),
\end{equation*}
with $\alpha_{1}=\paren{\begin{matrix} 
	1 & 0 \\
	-1 & 0
	\end{matrix}}+\paren{\begin{matrix} 
	0 & 0 \\
	1 & 0
	\end{matrix}}z
+\paren{\begin{matrix} 
	1 & 1 \\
	0 & 0
	\end{matrix}}\frac{z^2}{2}$,
$\alpha_{2}=\paren{\begin{matrix} 
	0 & 0 \\
	1 & 1
	\end{matrix}}+\paren{\begin{matrix} 
	0 & 0 \\
	-1 & 0
	\end{matrix}}z
+\paren{\begin{matrix} 
	-1 & -1 \\
	0 & 0
	\end{matrix}}\frac{z^2}{2}$,\\
$\alpha_{3}=\paren{\begin{matrix} 
	1 & 1 \\
	-1 & -1
	\end{matrix}}z
+\paren{\begin{matrix} 
	-1 & 1 \\
	0 & 0
	\end{matrix}}\frac{z^2}{2}$,
$\alpha_{4}=\paren{\begin{matrix} 
	0 & 0 \\
	1 & 0
	\end{matrix}}\frac{z^2}{2}$,
$\alpha_{5}=\paren{\begin{matrix} 
	0 & 0 \\
	0 & 1
	\end{matrix}}\frac{z^2}{2}$.
If  $E=span\llave{\alpha_{i}|1\leq i \leq 5}$ we obtain this step.

Step 3: The algebra $\Gamma$ is generated by $E$, $\complex\cdot \langle E\rangle=\Gamma$.

Since $\alpha_{1}+\alpha_{2}=I$ we have that $E=span\llave{I,\alpha_{1},\alpha_{3},\alpha_{4},\alpha_{5}}$. On the other hand,
{\small \begin{equation*}
	\alpha_{1}^2=\frac{4+2z^2+2z^3+z^4}{4}e_{11}+\frac{2z^2+z^4}{4}e_{12}+\frac{-2+2z-z^2+z^3}{2}e_{21}+\frac{z^3-z^2}{2}e_{22}
	\end{equation*}
	\begin{equation*}
	\alpha_{1}\alpha_{3}=-\frac{-4z+2z^2+z^4}{4}e_{11}+\frac{4z+2z^2+z^4}{4}e_{12}-\frac{2z-3z^2+z^3}{2}e_{21}+\frac{-2z+z^2+z^3}{2}e_{22}
	\end{equation*}
	\begin{equation*}
	\alpha_{1}\alpha_{4}=\frac{e_{11}z^4}{4}, \hspace{0.1cm}
	\alpha_{1}\alpha_{5}=\frac{e_{12}z^4}{4}, \hspace{0.1cm}
	\alpha_{3}\alpha_{1}=-\frac{z^4-4z^3}{4}e_{11}-\frac{z^4-2z^3}{4}e_{12}-\frac{2z^2+z^3}{2}e_{21}-\frac{z^3}{2}e_{22}
	\end{equation*}
	\begin{equation*}
	\alpha_{3}^2=\frac{z^4-6z^3}{4}e_{11}-\frac{2z^3+z^4}{4}e_{12}+\frac{z^3}{2}e_{21}-\frac{z^3}{2}e_{22}
	\end{equation*}
	\begin{equation*}
	\alpha_{3}\alpha_{4}=\frac{2z^3+z^4}{4}e_{11}-\frac{z^3}{2}e_{21},  \hspace{0.1cm}
	\alpha_{3}\alpha_{5}=\frac{2z^3+z^4}{4}e_{12}-\frac{z^3}{2}e_{22}
	\end{equation*}
	\begin{equation*}
	\alpha_{4}\alpha_{1}=\frac{2z^2+z^4}{4}e_{21}+\frac{z^4}{4}e_{22}, \hspace{0.1cm}
	\alpha_{4}\alpha_{3}=-\frac{z^4-2z^3}{4}e_{21}+\frac{2z^3+z^4}{4}e_{22}
	\end{equation*}
	\begin{equation*}
	\alpha_{4}^2=0, \hspace{0.1cm}
	\alpha_{4}\alpha_{5}=0, \hspace{0.1cm}
	\alpha_{5}\alpha_{1}=\frac{z^3-z^2}{2}e_{21}, \hspace{0.1cm}
	\alpha_{5}\alpha_{3}=-\frac{z^3}{2}e_{21}-\frac{z^3}{2}e_{22}, \hspace{0.1cm}
	\alpha_{5}\alpha_{4}=\frac{e_{21}z^4}{4}, \hspace{0.1cm}
	\alpha_{5}^2=\frac{e_{22}z^4}{4}.
	\end{equation*}}

Therefore, $e_{ij}z^3\in \complex\cdot \langle E\rangle$ and $e_{ij}z^4 \in \complex\cdot \langle E\rangle$ for $1\leq i,j \leq 2$ and using $\alpha_{4}, \alpha_{5}$ we obtain that 
$e_{ij}z^k \in \complex\cdot \langle E\rangle$ for $1\leq i,j \leq 2$ and $k\geq 3$. In particular $\complex\cdot \langle E\rangle =\Gamma$.

Step 4: The inclusion $\mathbb{A}\cap \oplus_{k=0}^2 M_{2}(\complex)[z]_{k}\subset E$.

Let $F\in \mathbb{A}\cap \oplus_{k=0}^2 M_{2}(\complex)[z]_{k}$ then there exists $\mathcal{L}=\mathcal{L}(x,\partial_{x})$ such that 
$\mathcal{L}\psi=\psi F$. We write 
\begin{equation*}
F(z)=\paren{\begin{matrix} 
	s_{0}^{11} & s_{0}^{12} \\
	s_{0}^{21} & s_{0}^{22}
	\end{matrix}}+
\paren{\begin{matrix} 
	s_{1}^{11} & s_{1}^{12} \\
	s_{1}^{21} & s_{1}^{22}
	\end{matrix}}z+
\paren{\begin{matrix} 
	s_{2}^{11} & s_{2}^{12} \\
	s_{2}^{21} & s_{2}^{22}
	\end{matrix}}z^2.
\end{equation*}	
After a computation we obtain that 
{\small	\begin{equation}\label{bispectral pair 3}
	\mathcal{L}=\left(\begin{matrix} 
	\frac{s_{0}^{11}x^{5}+( -s_{1}^{12}-4s_{0}^{11})x^{4}+( 3 s_{2}^{12}-4s_{2}^{11}+2s_{1}^{12}+4s_{0}^{11}) x^{3}+( -3s_{2}^{12}+12s_{2}^{11}+s_{1}^{21})x^2+( -s_{2}^{22}+2s_{2}^{21}-11s_{2}^{11}-2s_{1}^{21}) x-3s_{2}^{21}}{4x^{3}-4x^{4}+x^{5}}  \\ \\
	\frac{s_{0}^{21}x^4+(-s_{1}^{22}+s_{1}^{21}+s_{1}^{11}-4s_{0}^{21})x^3+(2s_{2}^{22}-4s_{2}^{21}-s_{2}^{12}+2s_{1}^{22}-4s_{1}^{21}
		-2s_{1}^{11}+4s_{0}^{21})x^2+(-s_{2}^{22}+11s_{2}^{21}+s_{2}^{11}+4s_{1}^{21})x-9s_{2}^{21}}{4x^2-4x^3+x^4}
	\end{matrix}\right.
	\end{equation}
	
	\begin{equation*}
	\left. \begin{matrix} 
	\frac{s_{0}^{12}x^6+(-s_{1}^{12}-4s_{0}^{12})x^5+(-s_{2}^{12}+4s_{1}^{12}+4s_{0}^{12})x^4+(s_{2}^{12}+s_{1}^{22}
		-4s_{1}^{12}-s_{1}^{11})x^3+(s_{2}^{22}-s_{2}^{12}+5s_{2}^{11}-2s_{1}^{22}+2s_{1}^{11})x^2
		+(-s_{2}^{22}-7s_{2}^{11})x-s_{2}^{21}}{4x^4-4x^5+x^6} \\ \\
	\frac{s_{0}^{22}x^5+(s_{1}^{12}-4s_{0}^{22})x^4+(-2s_{2}^{22}-s_{2}^{12}-2s_{1}^{12}+4s_{0}^{22})*x^3
		+(4s_{2}^{22}+3s_{2}^{12}-s_{1}^{21})x^2+(-3s_{2}^{22}+4s_{2}^{21}-s_{2}^{11}+2s_{1}^{21})x-5s_{2}^{21}}{4x^3-4x^4+x^5}
	\end{matrix}\right)
	\end{equation*}

	\begin{equation*}
	+	\left(\begin{matrix} 
	\frac{s_{1}^{11}x^3+(-s_{2}^{12}-2s_{1}^{11})x^2+s_{2}^{21}}{x^3-2x^2} &
	\frac{s_{1}^{12}x^3+(-s_{2}^{12}-2s_{1}^{12})x^2+2s_{2}^{12}x+s_{2}^{22}-s_{2}^{11}}{x^3-2x^2}\\ \\
	\frac{s_{1}^{21}x^2+(-s_{2}^{22}+s_{2}^{21}+s_{2}^{11}-2s_{1}^{21})x-2s_{2}^{21}}{x^2-2x} &
	\frac{s_{1}^{22}x^3+(s_{2}^{12}-2s_{1}^{22})x^2-s_{2}^{21}}{x^3-2x^2}
	\end{matrix}\right)\partial_{x}  +
	\left(\begin{matrix} 
	s_{2}^{11} & s_{2}^{12} \\ \\
	s_{2}^{21} & s_{2}^{22}
	\end{matrix}\right)\partial_{x}^{2}
	\end{equation*}}
and $s_{0}^{21}=s_{0}^{22}-s_{0}^{11}$, $s_{0}^{12}=0$, $s_{1}^{12}=s_{1}^{11}$, $s_{1}^{22}=-s_{1}^{11}$,
$s_{1}^{21}=s_{0}^{11}-s_{0}^{22}-s_{1}^{11}$, $s_{2}^{11}=\frac{s_{0}^{11}-s_{0}^{22}-s_{1}^{11}}{2}$,
$s_{2}^{12}=\frac{s_{1}^{12}+s_{0}^{11}-s_{0}^{22}}{2}$, then $F\in E$. 

Step 5: The inclusion $E\subset \mathbb{A}\cap\oplus_{k=0}^{2}M_{2}(\complex)[z]_{k}$

By the previous step we have Equation~\ref{bispectral pair 3} valid for every $F \in E$ and $(\mathcal{L} \psi )(x,z)=\psi(x,z)F(z)$, then $F \in  \mathbb{A}\cap\oplus_{k=0}^{2}M_{2}(\complex)[z]_{k}$.  

We are under the hypothesis of Theorem \ref{bis}, this concludes the proof of the assertion.\qed
\end{demos}

\section*{Conclusion and Final Comments}

In this article we characterized the bispectral triples 
associated to a certain class of matrix-valued eigenfuctions.
Furthermore, we established important  properties of the full rank $1$
algebras  as a model of some bispectral algebras. These properties include the fact that they are Noetherian and finitely generated. An important role was played by the Ad-condition due to the fact that the matrix-valued operators were acting from opposite directions. 

Natural questions arise when searching for a characterization of the bispectral partners 
to a given operator. Indeed, there is a family of maps $\llave{P_k}_{k\in \natu}$ with the  translation and product properties of Theorem \ref{Product Pk} and Lemma \ref{translation} which generate the algebra in  Example \ref{first}. Is it possible to do this for the general case? Or, would this be possible at least for  the bispectral partners of a given Schr\"{o}dinger operator? 


\section*{Acknowledgments}
The authors acknowledge the financial support provided by the  {\em Coordena\c{c}\~{a}o de Aperfei\c{c}oamento de Pessoal de N\'{\i}vel Superior} (CAPES),  {\em Conselho Nacional de Desenvolvimento Cient\'{\i}fico e Tecnol\'{o}gico} (CNPq), and {\em Funda\c{c}\~{a}o Carlos Chagas Filho de Amparo \`{a} Pesquisa do Estado do Rio de Janeiro} (FAPERJ).

BVDC was supported by CAPES grants 88882 332418/2019-01 as well as IMPA.
JPZ was supported by CNPq grants 302161 and 47408, as well as by FAPERJ under the program \textit{Cientistas do Nosso Estado}.

\appendix
\section{The Ad-Condition and Polynomial Eigenvalues}\label{adcond}
We close this article with a generalization of a key lemma from the work of \cite{duistermaat1986differential} to the noncommutative case.
\begin{prop}\label{degad}
	Let $\mathcal{L}=\mathcal{L}(x,\partial_{x})=\sum_{j=0}^{l}L_{j}(x)\partial_{x}^{j}$, $\Theta=\Theta(x,\partial_{x})=\sum_{s=0}^{m}\theta_{s}(x)\partial_{x}^{s}$
	then  $\corch{L_{l},\theta_{m}}=0$ implies  $\deg_{\partial_{x}}( ad (\mathcal{L})(\Theta))\leq m+l-1$ and $\corch{L_{l},\theta_{m}}\neq 0$
	implies  $\deg_{\partial_{x}}( ad (\mathcal{L})(\Theta))= m+l$ .	
\end{prop}
\begin{demos}
	By definition $(ad \hspace{0.1 cm} \mathcal{L})(\Theta)= \corch{\mathcal{L},\Theta}= \mathcal{L}\Theta -\Theta \mathcal{L}$ then 
	\begin{equation*}
	(ad \mathcal{L})(\Theta)= \mathcal{L}\Theta-\Theta \mathcal{L}
	= \paren{\sum_{j=0}^{l}L_{j}(x)\partial_{x}^{j}}\paren{\sum_{s=0}^{m}\theta_{s}(x)\partial_{x}^{s}}-
	\paren{\sum_{s=0}^{m}\theta_{s}(x)\partial_{x}^{s}}\paren{\sum_{j=0}^{l}L_{j}(x)\partial_{x}^{j}}
	\end{equation*}
	\begin{equation*}
	= \sum_{j=0}^{l}\sum_{s=0}^{m}L_{j}\partial_{x}^{j}(\theta_{s}\partial_{x}^{s})-\sum_{s=0}^{m}\sum_{j=0}^{l}
	\theta_{s}\partial_{x}^{s}(L_{j}\partial_{x}^{j})
	=  \sum_{j=0}^{l}\sum_{s=0}^{m}L_{j}\sum_{k=0}^{j}\paren{	\begin{matrix} 
		j \\
		k  
		\end{matrix}} \theta_{s}^{(k)}\partial_{x}^{j-k+s}-\sum_{s=0}^{m}\sum_{j=0}^{l}
	\theta_{s}\sum_{r=0}^{s}\paren{	\begin{matrix} 
		s \\
		r  
		\end{matrix}}L_{j}^{(r)}\partial_{x}^{s-r+j}
	\end{equation*}	
	\begin{equation*}
	=  \sum_{k=0}^{l}\sum_{s=0}^{m}\sum_{j=k}^{l}L_{j}\paren{	\begin{matrix} 
		j \\
		k  
		\end{matrix}} \theta_{s}^{(j-k)}\partial_{x}^{k+s}-\sum_{s=0}^{m}\sum_{j=0}^{l}
	\sum_{j=s}^{m}\theta_{j}\paren{	\begin{matrix} 
		j \\
		s  
		\end{matrix}}L_{k}^{(j-s)}\partial_{x}^{k+s}
	\end{equation*}
	\begin{equation*}
	=  \sum_{k=0}^{l}\sum_{s=0}^{m}\paren{\sum_{j=k}^{l}L_{j}\paren{	\begin{matrix} 
			j \\
			k  
			\end{matrix}} \theta_{s}^{(j-k)}-
		\sum_{j=s}^{m}\theta_{j}\paren{	\begin{matrix} 
			j \\
			s  
			\end{matrix}}L_{k}^{(j-s)}}\partial_{x}^{k+s}
	\end{equation*}
	\begin{equation*}
	=\sum_{r=0}^{m+l}\paren{\sum_{k+s=r, 0\leq k \leq l, 0\leq s\leq m}\paren{\sum_{j=k}^{l}L_{j} \paren{	\begin{matrix} 
				j \\
				k  
				\end{matrix}} \theta_{s}^{(j-k)}-
			\sum_{j=s}^{m}\theta_{j}\paren{	\begin{matrix} 
				j \\
				s  
				\end{matrix}}L_{k}^{(j-k+r)}}}\partial_{x}^{r}
	\end{equation*}
	\begin{equation*}
	=\sum_{r=0}^{m+l}a_{r}\partial_{x}^{r}
	\end{equation*}
	with \begin{equation*}
	a_{r}=\sum_{k+s=r, 0\leq k \leq l, 0\leq s\leq m}\paren{\sum_{j=k}^{l}L_{j} \paren{	\begin{matrix} 
			j \\
			k  
			\end{matrix}} \theta_{s}^{(j-k)}-
		\sum_{j=s}^{m}\theta_{j}\paren{	\begin{matrix} 
			j \\
			s  
			\end{matrix}}L_{k}^{(j-k+r)}}
	\end{equation*} 
	in particular $a_{m+l}=L_{l}\theta_{m}-\theta_{m}L_{l}=\corch{L_{k},\theta_{m}}$, hence 
	$\corch{L_{k},\theta_{m}}=0$ implies  $\deg_{\partial_{x}}( ad (\mathcal{L})(\theta))\leq m+l-1$ and $\corch{L_{k},\theta_{m}}\neq 0$
	implies  $\deg_{\partial_{x}}( ad (\mathcal{L})(\theta))= m+l$ .\qed
\end{demos}

The previous proposition implies that if $k+s=m+l-1$, \hspace{0.1cm} $0\leq k\leq l$, \hspace{0.1cm} $0\leq s \leq m$ then $(k,s)=(l-1,m)$ or $(k,s)=(l,m-1)$, therefore
\begin{equation*}
a_{m+l-1}=\sum_{j=l-1}^{l}L_{j}\paren{	\begin{matrix} 
	j \\
	j-1  
	\end{matrix}}\theta_{m}^{(j-l+1)}-\theta_{m}L_{l-1}+L_{l}\theta_{m-1}-\sum_{j=m-1}^{m}\theta_{j}\paren{	\begin{matrix} 
	j \\
	m-1  
	\end{matrix}}L_{l}^{(j-m+1)}  
\end{equation*}
\begin{equation*}
=L_{l-1}\theta_{m}+lL_{l}\theta_{m}^{'}-\theta_{m}L_{l-1}+L_{l}\theta_{m-1}-\theta_{m-1}L_{l}
-m\theta_{m}L_{l}^{'}=\corch{L_{l-1},\theta_{m}}+\corch{L_{l},\theta_{m-1}}+lL_{l}\theta_{m}^{'}-m\theta_{m}L_{l}^{'}.
\end{equation*}
If $\corch{L_{l-1},\theta_{m}}=0$, $\corch{L_{l},\theta_{m-1}}=0$, then 
\begin{equation*}
a_{m+l-1}=lL_{l}\theta_{m}^{'}-m\theta_{m}L_{l}^{'}.
\end{equation*}
In particular if $m=0$, $\Theta=\theta_{0}$ and $\corch{L_{l-1},\theta_{0}}=0$, then 
\begin{equation*}
a_{l-1}=lL_{l}\theta_{0}^{'}.
\end{equation*}
If we assume the system of equations \ref{bispec} we obtain:
\begin{equation*}
(ad \hspace{0.1 cm}\mathcal{L})(\theta) \psi= \corch{\mathcal{L},\theta}\psi=\paren{\mathcal{L}\theta-\theta \mathcal{L}}\psi =\mathcal{L}(\theta \psi)-\theta \mathcal{L}\psi = \mathcal{L}(\psi \mathcal{B})- \theta  \psi F
\end{equation*}
\begin{equation*}
=(\mathcal{L}\psi)\mathcal{B}- (\psi \mathcal{B})F=(  \psi F)\mathcal{B}-\psi \mathcal{B} F =\psi \corch{F,\mathcal{B}}=\psi (ad \hspace{0.1 cm} F)(\mathcal{B}).
\end{equation*}

Now we prove by induction that 
\begin{equation*}
(ad \mathcal{L})^{r}(\theta)\psi=\psi (ad \hspace{0.1 cm} F)^{r}(\mathcal{B}),
\end{equation*}
for all $r\in \ente_{+}$.

The claim is clear for $r=1$. Assume the condition for $r$ and consider the case $r+1$, then 
\begin{equation*}
\psi (ad \hspace{0.1 cm}F)^{r+1}(\mathcal{B})=\psi (ad \hspace{0.1 cm} F)(ad \hspace{0.1 cm} F)^r (\mathcal{B})= (F \psi) (ad \hspace{0.1 cm} F)^r (\mathcal{B})-F (\psi) (ad \hspace{0.1 cm} F)^r (\mathcal{B}) 
\end{equation*}
\begin{equation*}
=(F \psi)(ad \hspace{0.1 cm} F)^r (\mathcal{B})-F (ad \mathcal{L})^r (\theta)\psi= (\mathcal{L}\psi)(ad \hspace{0.1 cm} F)^r(\mathcal{B})-F (ad \mathcal{B})^{r}(\theta)\psi
\end{equation*}
\begin{equation*}
=\mathcal{L}(\psi (ad \hspace{0.1 cm} F)^r (\mathcal{B}))-(ad \mathcal{L})^r(\theta)(F \psi)= \mathcal{L} (ad \mathcal{L})^r(\theta)\psi - (ad \mathcal{L})^{r}(\theta)(\mathcal{L}\psi)
\end{equation*}
\begin{equation*}
=(\mathcal{L}(ad \mathcal{L})^r(\theta)-(ad \mathcal{L})^r (\theta)\mathcal{L})\psi= (ad \mathcal{L})(ad \mathcal{L})^{r}(\theta)\psi= (ad \mathcal{\mathcal{L}})^{r+1}(\theta)\psi.
\end{equation*}

If $\deg_{\partial_{z}}\mathcal{B}=m$ and $F$ is scalar we use the Proposition \ref{degad} to conclude $(ad \hspace{0.1 cm} \mathcal{L})^{m+1}(\theta)\psi=\psi (ad \hspace{0.1 cm} F)^{m+1}(\mathcal{B})=0$, similarly if $\deg \mathcal{L}=l$ then  $\deg_{\partial_{x}} (ad \mathcal{L})^{m+1}(\theta)\leq (m+1)(l-1)$, in our case $\deg_{\partial_{x}} (ad \mathcal{L})^{m+1}(\theta)\leq m+1< \infty$
since $\psi(\cdot,z)\in \ker((ad \mathcal{L})^{m+1}(\theta))$ for every $z\in \complex$ and $\llave{\psi(\cdot,z)}_{z\in \complex}$ is a linearly independent set and $\dim \ker((ad \mathcal{L})^{m+1}(\theta))\leq \deg_{\partial_{x}} (ad \mathcal{L})^{m+1}(\theta)$ if $(ad \mathcal{L})^{m+1}(\theta)\neq 0$ we have that $(ad \mathcal{L})^{m+1}(\theta)=0$.

Finally we claim that if $\mathcal{L}=\sum_{j=0}^{l} L_{j}\partial_{x}^j$ with $L_{l}\in \complex \setminus \llave{0}$  and $L_{l-1}=0$, then $\mathbf{coeff}((ad \mathcal{L})^{k+1}(\theta),\partial_{x},(k+1)(l-1))=
(lL_{l})^{k+1}\theta^{(k+1)}(0)$ for every $k\in \natu$.

The claim is obvious for $k=0$. If we assume that the claim is valid for $k$, then 
\begin{equation*}
\mathbf{coeff}((ad \mathcal{L})^{k+2}(\theta),\partial_{x},(k+2)(l-1))=\mathbf{coeff}((ad \mathcal{L})(ad \mathcal{L})^{k+1}(\theta),\partial_{x},(k+2)(l-1))
\end{equation*}
\begin{equation*}
=lL_{l}\partial_{x}((lL_{l})^{k+1}\theta^{(k+1)})=(lL_{l})(lL_{l})^{k+1}\theta^{(k+2)}=(lL_{l})^{k+2}\theta^{(k+2)}
\end{equation*} because $L_{l}$ is constant and scalar. 

Since  $(ad \mathcal{L})^{m+1}(\theta)=0$ we have that $\theta^{(m+1)}=0$ and $\theta$ has to be a polynomial with $\deg \theta\leq m$.
{\small
\bibliographystyle{alpha}
\bibliography{new2019}}

\begin{thebibliography}{ABvW17}

\bibitem[ABvW17]{Pidecom}
PN~Anh, GF~Birkenmeier, and L~van Wyk.
\newblock Peirce decompositions, idempotents and rings.
\newblock {\em arXiv preprint arXiv:1702.05261}, 2017.

\bibitem[BGK08]{BGK09}
Maarten Bergvelt, Michael Gekhtman, and Alex Kasman.
\newblock Spin calogero particles and bispectral solutions of the matrix {KP}
  hierarchy.
\newblock {\em Mathematical Physics Analysis and Geometry}, 12, 07 2008.

\bibitem[BL08]{BL08}
Carina Boyallian and Jose Liberati.
\newblock Matrix-valued bispectral operators and quasideterminants.
\newblock {\em Journal of Physics A: Mathematical and Theoretical}, 41:365209,
  08 2008.

\bibitem[Cas17]{Casper2017}
William~R. Casper.
\newblock {\em Bispectral Operator Algebras}.
\newblock University of Washington, 2017.
\newblock Thesis (Ph.D.)--University of Washington, 2017-06.

\bibitem[DG86]{duistermaat1986differential}
J.J. Duistermaat and F.A. Gr{\"{u}}nbaum.
\newblock Differential equations in the spectral parameter.
\newblock {\em Communications in Mathematical Physics}, 103(2):177--240, 1986.

\bibitem[DG04]{DuranGrunbaum2004}
A.~J. Dur{\'{a}}n and F.A. Gr{\"{u}}nbaum.
\newblock Orthogonal matrix polynomials satisfying second-order differential
  equations.
\newblock {\em Int. Math. Res. Not.}, 2004(10), 2004.

\bibitem[GHY15]{GHY17}
Joel Geiger, E.~Horozov, and Milen Yakimov.
\newblock Noncommutative bispectral darboux transformations.
\newblock {\em Transactions of the American Mathematical Society}, 369, 08
  2015.

\bibitem[GI03]{GI03}
Alberto Gr{\"{u}}nbaum and Plamen Iliev.
\newblock A noncommutative version of the bispectral problem.
\newblock {\em J. Comput. Appl. Math.}, 161:99--118, 12 2003.

\bibitem[Gr{\"{u}}14]{Grfrm-4}
F.~Alberto Gr{\"{u}}nbaum.
\newblock Some noncommutative matrix algebras arising in the bispectral
  problem.
\newblock {\em SIGMA Symmetry Integrability Geom. Methods Appl.}, 10:078, 2014.

\bibitem[Hor02]{horozov2002bispectral}
E.~Horozov.
\newblock Bispectral operators of prime order.
\newblock {\em Communications in Mathematical Physics}, 231(2):287--308, 2002.

\bibitem[Ili99]{iliev1999discrete}
P.~Iliev.
\newblock {\em Discrete versions of the Kadomtsev-Petviashvili hierarchy and
  the bispectral problem}.
\newblock PhD thesis, Dissertation, 1999.

\bibitem[Kas15]{Kas15}
Alex Kasman.
\newblock Bispectrality of $n$-component {KP} wave functions: A study in
  non-commutativity.
\newblock {\em Symmetry, Integrability and Geometry: Methods and Applications},
  11, 05 2015.

\bibitem[Sta85]{Stafford}
J.~T. Stafford.
\newblock On the ideals of a noetherian ring.
\newblock {\em Transactions of the American Mathematical Society},
  289(1):381--392, 1985.

\bibitem[SZ01]{SZ01}
Alexander Sakhnovich and Jorge Zubelli.
\newblock Bundle bispectrality for matrix differential equations.
\newblock {\em Integral Equations Operator Theory}, 41(4):472, 2001.

\bibitem[Wil93]{wilson}
George Wilson.
\newblock Bispectral commutative ordinary differential operators.
\newblock {\em J. Reine Angew. Math.}, 442:177--204, 1993.

\bibitem[ZM91]{zubelli1991differential}
J.P. Zubelli and F.~Magri.
\newblock Differential equations in the spectral parameter, {Darboux}
  transformations and a hierarchy of master symmetries for {KdV}.
\newblock {\em Communications in Mathematical Physics}, 141(2):329--351, 1991.

\bibitem[Zub90]{zubelli1990differential}
J.P. Zubelli.
\newblock Differential equations in the spectral parameter for matrix
  differential operators.
\newblock {\em Physica D: Nonlinear Phenomena}, 43(2-3):269--287, 1990.

\bibitem[Zub92a]{Zubelli1992a}
Jorge~P. Zubelli.
\newblock On a zero curvature problem related to the {Z}{S}-{A}{K}{N}{S}
  operator.
\newblock {\em J. Math. Phys.}, 33(11):3666--3675, 1992.

\bibitem[Zub92b]{Zubelli1992}
Jorge~P. Zubelli.
\newblock On the polynomial $\tau$-functions for the {K}{P} hierarchy and the
  bispectral property.
\newblock {\em Lett. Math. Phys.}, 24(1):41--48, 1992.

\bibitem[Zub92c]{Zubelli1992b}
Jorge~P. Zubelli.
\newblock Rational solutions of nonlinear evolution equations, vertex
  operators, and bispectrality.
\newblock {\em J. Differential Equations}, 97(1):71--98, 1992.

\bibitem[ZVS00]{Zubelli2000}
Jorge~P. Zubelli and D.~S. Valerio~Silva.
\newblock Rational solutions of the master symmetries of the {K}d{V} equation.
\newblock {\em Communications in Mathematical Physics}, 211(1):85--109, 2000.

\end{thebibliography}
\end{document}